\documentclass[11pt, reqno]{amsart}
\usepackage{lmodern}
\usepackage{amsmath, amsthm, amssymb, amsfonts}
\usepackage[normalem]{ulem}
\usepackage{hyperref}
\usepackage[all,cmtip]{xy}
\usepackage{verbatim}
\usepackage{nccmath}
\usepackage{stmaryrd}
\usepackage{caption}
\usepackage{mathrsfs}
\usepackage{mathtools}
\usepackage{esvect}
\usepackage{cite}
\usepackage{bbm}
\usepackage{eucal}

\usepackage{mathbbol}
\usepackage{tabularx}
\usepackage[toc,page]{appendix}

\usepackage{tikz-cd}

\theoremstyle{plain}
\newtheorem{thm}{Theorem}[section]
\newtheorem{cor}[thm]{Corollary}

\newtheorem{lemma}[thm]{Lemma}
\newtheorem{prop}[thm]{Proposition}

\newtheorem*{thmn}{Theorem}
\newtheorem*{defnn}{Definition}

\theoremstyle{definition}
\newtheorem{defn}[thm]{Definition}

\makeatletter
\newcommand\ackname{Acknowledgements}
\if@titlepage
\newenvironment{acknowledgements}{%
	\titlepage
	\null\vfil
	\@beginparpenalty\@lowpenalty
	\begin{center}%
		\bfseries \ackname
		\@endparpenalty\@M
\end{center}}%
{\par\vfil\null\endtitlepage}
\else

\fi
\makeatother

\theoremstyle{remark}
\newtheorem{rmk}[thm]{Remark}

\newcommand{\BC}{{\mathbb{C}}}

\newcommand{\BE}{{\mathbb{E}}}

\newcommand{\BL}{{\mathbb{L}}}
\newcommand{\BM}{{\mathbb{M}}}

\newcommand{\BQ}{{\mathbb{Q}}}
\newcommand{\BR}{{\mathbb{R}}}

\newcommand{\BT}{{\mathbb{T}}}

\newcommand{\BZ}{{\mathbb{Z}}}

\newcommand{\CA}{{\mathcal A}}

\newcommand{\CC}{{\mathcal C}}
\newcommand{\CD}{{\mathcal D}}
\newcommand{\CE}{{\mathcal E}}

\newcommand{\CI}{{\mathcal I}}

\newcommand{\CK}{{\mathcal K}}
\newcommand{\CL}{{\mathcal L}}
\newcommand{\CM}{{\mathcal M}}

\newcommand{\CO}{{\mathcal O}}
\newcommand{\CP}{{\mathcal P}}

\newcommand{\FM}{{\mathfrak{M}}}

\newcommand{\ch}{{\mathrm{ch}}}

\DeclareFontFamily{OT1}{rsfs}{}
\DeclareFontShape{OT1}{rsfs}{n}{it}{<-> rsfs10}{}
\DeclareMathAlphabet{\curly}{OT1}{rsfs}{n}{it}

\newcommand{\p}{\mathbb{P}}

\newcommand\Spec{\operatorname{Spec}}

\newcommand{\Mbar}{{\overline M}}

\newcommand{\Pic}{\mathop{\rm Pic}\nolimits}

\newcommand{\ev}{{\mathrm{ev}}}
\newcommand{\evsf}{{\mathsf{ev}}}

\newcommand{\br}{\mathsf{br}}
\newcommand{\Aut}{\mathrm{Aut}}

\newcommand{\FC}{\mathfrak{C}}

\begin{document}
	
	\title[Gromov--Witten/Hurwitz wall-crossing]
	{Gromov--Witten/Hurwitz wall-crossing}

\author{Denis Nesterov}
\address{ETH Z\"urich, Departement Mathematik}
\email{denis.nesterov@math.ethz.ch}
	\maketitle
	\begin{abstract}
		

		For a target variety $X$ and a nodal curve $C$, we introduce a one-parameter stability condition, termed $\epsilon$-admissibility, for maps from nodal curves to  $X\times C$. If $X$ is a point, $\epsilon$-admissibility interpolates between moduli spaces of  stable maps to $C$ relative to some fixed points and moduli spaces of admissible covers with arbitrary ramifications over the same fixed points and simple ramifications elsewhere on $C$.

		Using Zhou's entangled tails, we prove wall-crossing formulas relating  invariants for different values of $\epsilon$. If $X$ is a surface, we use this wall-crossing in conjunction with author's quasimap wall-crossing to show that the relative Pandharipande--Thomas/Gromov--Witten correspondence of $X\times C$ and Ruan's extended crepant resolution conjecture of the pair $X^{[n]}$ and $[X^{(n)}]$ are equivalent up to explicit wall-crossings. We thereby prove the crepant resolution conjecture for 3-point genus-0 invariants in all classes, if $X$ is a toric del Pezzo surface.

\end{abstract}
\setcounter{tocdepth}{1}
\tableofcontents
\section{Introduction}
\subsection{Overview}
Inspired by the theory of quasimaps to GIT quotients of \cite{CFKM}, a theory of quasimaps to moduli spaces of sheaves was introduced in \cite{N}. When applied to Hilbert schemes of $n$-points $S^{[n]}$ of a surface $S$, moduli spaces of $\epsilon$-stable quasimaps interpolate between moduli spaces of stable maps to  $S^{[n]}$  and Hilbert schemes of 1-dimensional subschemes of the relative geometry $S\times C_{g,N}/\Mbar_{g,N}$, where $C_{g,N} \rightarrow \Mbar_{g,N}$ is the universal curve of a moduli space of stable marked curves, 
\begin{equation} \label{qmwall}
	\xymatrix{
		\Mbar_{g,N}(S^{[n]},\beta)  \ar@{<-->}[r]|-{\epsilon } &\mathrm{Hilb}_{n,\check{\beta}}(S\times C_{g,N}).}
\end{equation}

This interpolation gives rise to wall-crossing formulas, which therefore relate Gromov--Witten (GW) theory of $S^{[n]}$ and relative Donaldson--Thomas (DT) theory of $S\times C_{g,N}/\Mbar_{g,N}$. Alongside with results of \cite{NK3}, the quasimap wal-crossing was used to prove various correspondences, among which is the wall-crossing part of Igusa cusp form conjecture of \cite{OPa}. For more details, we refer to \cite{N,NK3}.

In this article we introduce a notion of $\epsilon$-admissibility, depending on a parameter $\epsilon \in (0,1] \subset \BR_{>0}$, for maps 
$$ P \rightarrow X\times (C,\mathbf{x})$$ 
relative to $X\times \mathbf x$, where $P$ is a nodal curve, $(C, \mathbf x)$ is a marked nodal curve and $X$ is a smooth projective variety. 

As the value of $\epsilon$ varies, moduli spaces of $\epsilon$-admissible maps interpolate between moduli spaces of stable twisted maps to an orbifold symmetric product $[X^{(n)}]$ of \cite{AGV} and moduli spaces of stable maps with possibly disconnected sources to the relative geometry $X\times C_{g,N}/\Mbar_{g,N}$,  
\[\xymatrix{
	\CK_{g,N}([X^{(n)}],\beta) \ar@{<-->}[r]|-{\epsilon } &\Mbar_{\mathsf h}(X\times C_{g,N}, (\gamma,n)),}
\]
such that the various discrete data on both sides, like the genus of a curve or the degree of a map, determine each other, as is explained in Section \ref{Relation1}.  Note that we consider the notion of extended degree for maps in $\CK_{g,N}([X^{(n)}],\beta)$, introduced in \cite{BG};  after passing to coarse curves, it permits simple ramifications away from marked points. We also allow rational bridges for target curves in $\Mbar_{\mathsf h}(X\times C_{g,N}, (\gamma,n))$ without explicitly indicating this in the notation of the spaces.\footnote{The notation also does not indicate that we work relatively to $\Mbar_{g,N}$ for the sake of brevity.}

Using Zhou’s theory of entangled tails from \cite{YZ}, we establish wall-crossing formulas which relates the associated invariants for different values of  $\epsilon \in (0,1]$. This wall-crossing is completely analogous to the quasimap wall-crossing. The result is an equivalence of the orbifold GW theory of $[X^{(n)}]$ and the GW theory of the relative geometry $X\times C_{g,N}/\Mbar_{g,N}$ with relative insertions for an arbitrary smooth projective target $X$, which can be expressed in terms of a change of variables applied to certain generating series. The change of variables involves so-called \textit{I-functions}, which are defined via the localised GW theory of $X\times \p^1$ with respect to the $\BC^*$-action coming from the $\p^1$-factor.
This wall-crossing can be called a Gromov--Witten/Hurwitz (\textsf{GW/H}) wall-crossing, because if $X$ is a point, the moduli spaces of $\epsilon$-admissible maps interpolates between Gromov--Witten and Hurwitz spaces of a curve $C$. 
\\

\begin{figure} [h!]
	\scriptsize
	\[
	\begin{picture}(200,75)(-30,-50)
		\thicklines
		\put(25,-25){\line(1,0){30}}
		\put(95,-25){\line(1,0){30}}
		\put(25,-25){\line(0,1){40}}
		\put(25,30){\makebox(0,0){\textsf{Quasimap}}}
		\put(25,20){\makebox(0,0){\textsf{wall-crossing}}}
		\put(25,35){\line(0,1){40}}
		\put(125,-25){\line(0,1){40}}
		\put(125,30){\makebox(0,0){\textsf{GW/H}}}
		\put(125,20){\makebox(0,0){\textsf{wall-crossing}}}
		\put(25,75){\line(1,0){30}}
		\put(95,75){\line(1,0){30}}
		\put(125,35){\line(0,1){40}}
		\put(140,85){\makebox(0,0){$\mathsf{GW}_{\mathrm{orb}}([S^{(n)}])$}}
		\put(15,85){\makebox(0,0){$\mathsf{GW}(S^{[n]})$}}
		\put(75,75){\makebox(0,0){\textsf{CRC}}}
		\put(5,-35){\makebox(0,0){$\mathsf{PT_{rel}}(S\times C_{g,N})$}}
		\put(150,-35){\makebox(0,0){$\mathsf{GW_{rel}}(S\times C_{g,N})$}}
		\put(75,-25){\makebox(0,0){\textsf{PT/GW}}}
	\end{picture}
	\]
	\caption{The Square}
	\label{square}
\end{figure}
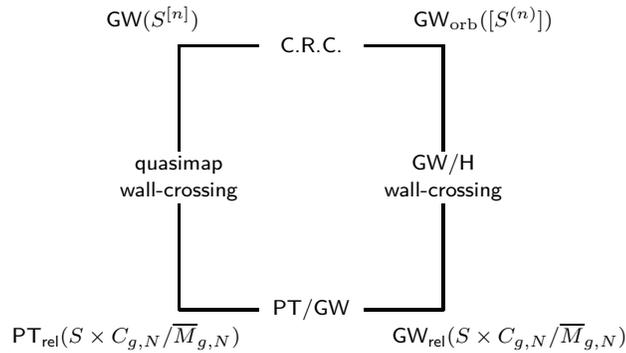

In conjunction with the quasimap wall-crossing of \cite{N}, \textsf{GW/H} wall-crossing establishes the square of theories for a smooth surface $S$, illustrated in Figure \ref{square}. The square relates the crepant resolution conjecture (\textsf{CRC}), proposed in \cite{YR} and refined in \cite{BG,CCIT}, and Pandharipande-- Thomas/Gromov--Witten correspondence (\textsf{PT/GW}),
 proposed in \cite{MNOP1,MNOP1}. The square has some similarities with Landau--Ginzburg/Calabi--Yau correspondence, as it is explained in Section \ref{LG/CY}.

With the help of the square, we establish the following results: 
	\begin{itemize}
	
	\item the 3-point genus-0 crepant resolution conjecture in the sense of \cite{BG} for the pair $S^{[n]}$ and $[S^{(n)}]$ in all classes, if $S$ is a toric del Pezzo surface. 
	\item a geometric origin of $y=-e^{iu}$ in \textsf{PT/GW} through $\mathsf{CRC}$.
\end{itemize}
Moreover, a cycle-valued version of the wall-crossing should have applications in the theory of double ramifications cycles of \cite{JPPZ}, comparison results for the TQFT's from \cite{Ca07} and \cite{BP}, etc.  This will be addressed in a future work. 

Various instances of the vertical sides of the square were studied on the level of invariants in numerous articles, mainly for $\BC^2$ and $\CA_m$ - \cite{OP10}, \cite{OP10b}, \cite{OP10c}, \cite{BP}, \cite{PT19}, \cite{Mau}, \cite{MO}, \cite{Che09}, etc. The wall-crossings provide a geometric justification for these phenomena. 

	\vspace{1 in}

\subsection{Analogy}
Let us  illustrate how the theory of quasimaps sheds light on a seemingly unrelated theme of admissible covers. 
\subsubsection{Moduli spaces of $\epsilon$-stable quasimaps}
A map from a nodal curve $C$, \[f\colon C \rightarrow S^{[n]},\]
is determined by its graph \[\Gamma_{f} \subset S\times C,\]
which is associated to $f$ via the construction of $S^{[n]}$ as a moduli space of $0$-dimensional subschemes on $S$. 
If the curve $C$ varies, the pair $(C, \Gamma_f)$ can degenerate in two ways: 
\begin{itemize}
	\item[(i)] the curve $C$ degenerates,
	\item[(ii)] the graph $\Gamma_f$ degenerates. 
\end{itemize}  
By a degeneration of $\Gamma_f$ we mean that it becomes non-flat\footnote{A 1-dimensional subscheme $\Gamma \subset S\times C$ is a graph, if and only if it is flat.} over $C$ as a subscheme of $S\times C$, which is due to
\begin{itemize}
	\item floating points,
	\item non-dominant components. 
\end{itemize}
Two types of degenerations of a pair $(C, \Gamma_f)$ are related. GW theory of $S^{[n]}$ proposes that $C$ sprouts out a rational tail ($C$ degenerates), whenever non-flatness arises ($\Gamma_f$ degenerates). DT theory, on the other hand, allows non-flatness, since it is interested in arbitrary 1-dimensional subschemes, thereby restricting degenerations of $C$ to semistable ones (no rational tails).

A non-flat graph $\Gamma$ does not define a map to $S^{[n]}$, but it defines a quasimap to $S^{[n]}$. Hence the motto of quasimaps:
\[\textit{Trade rational tails for non-flat points and vice versa.}\]
The idea of $\epsilon$-stability is to allow both rational tails and non-flat points, restricting their degrees.
The moduli spaces involved in (\ref{qmwall}) are given by the extremal values of $\epsilon$. 
\subsubsection{Moduli spaces of $\epsilon$-admissible maps}
The motto of the $\mathsf{GW/H}$ wall-crossing is the following one:
\[\textit{Trade rational tails for branching points and vice versa.}\] 
Let us explain what we mean by making an analogy with quasimaps. Let 
\[f\colon P \rightarrow C\] 
be an admissible cover with simple ramifications introduced in \cite[Chapter 4]{HM82}. If the curve $C$ varies, the pair $(C,f)$ can degenerate in two ways: 
\begin{itemize}
\item[(i)] the curve $C$ degenerates,
\item[(ii)] the cover $f$ degenerates.
\end{itemize}
The degenerations of $f$ arise due to
\begin{itemize} 
\item ramifications of higher order,
\item contracted components and singular points mapping to smooth locus.
\end{itemize}
As previously, these two types of degenerations of a pair $(C,f)$ are related. Hurwitz theory of a varying curve $C$ proposes that $C$ sprouts out rational tails, whenever $f$ degenerates in the sense above. On the other hand, GW theory of a varying curve $C$ allows $f$ to degenerate and therefore restricts the degenerations of $C$ to semistable ones. The purpose of $\epsilon$-admissible maps is to interpolate between these Hurwitz and Gromov--Witten cases. 
\subsection{Definition}
We will now outline the definition of $\epsilon$-admissibility. Let $f\colon P\rightarrow  C$ be a degree $n$ map between nodal curves, such that it is admissible at nodes (see \cite[Chapter 4]{HM82} for the admissibility) and $g(P)=\mathsf h$, $g(C)=g$. We allow $P$ to be disconnected, requiring that each connected component is mapped non-trivially to $C$. Following \cite{FP}, we define the branching divisor 
\[\br(f) \in \mathrm{Div}(C),\] 
it is an effective divisor which measures the degree of ramification away from nodes and the genera of contracted components of $P$. If $C$ is smooth, the branching divisor $\br(f)$ is given by associating to the 0-dimensional complex 
\[Rf_*[f^*\Omega_C \rightarrow \Omega_P ]\]
its support weighted by Euler characteristics. Otherwise, we need to take
the part of the support which is contained in the regular locus of $C$. More intrinsically, we can use logarithmic cotangent bundles to achieve this.

Using the branching divisor $\br$, we  define $\epsilon$-admissibility by the weighted stability of the pair $(C,\br(f))$, considered in \cite{Ha03}. A similar stability was studied in \cite{D}, where the source curve $P$ is allowed to have more degenerate singularities instead of contracted components. However, the moduli spaces of  \cite{D} do not have a perfect obstruction theory.

\begin{defnn}  Let $\epsilon \in (0,1]$. A map $f$ is $\epsilon$-admissible, if 
\begin{itemize}
	\item[(i)] $\omega_C(\epsilon \cdot \br(f))$ is ample,
	\item [(ii)] $\forall p \in C$,  $\mathrm{mult}_p(\br(f))\leq 1/\epsilon$, 
	\item[(iii)] the group of automorphisms of $f$ is finite. 
\end{itemize}
\end{defnn}


One can readily verify that for $\epsilon=1$, an $\epsilon$-admissible map is an admissible cover with simple ramifications. For a value of $\epsilon$ such that $\epsilon\ll 1$, which we denote by $0^+$, an $\epsilon$-admissible map is a stable\footnote{When the target curve $C$ is singular, we assume that all maps are admissible over nodes.} map, such that the target curve $C$ is semistable. Hence $\epsilon$-admissibility provides an interpolation between the moduli space of admissible covers with simple ramifications $Adm_{g,\mathsf h, n}$ and the moduli space of stable maps $\Mbar_{\mathsf h}(C_g,n)$ of $C_g/\Mbar_{g}$,
\[\xymatrix{
A dm_{g,\mathsf h, n} \ar@{<-->}[r]|-{\epsilon } &\Mbar_{\mathsf h}(C_g,n).}
\]
After introducing markings $\mathbf x=(x_1,\dots, x_N)$ on $C$ and requiring maps to be admissible over these markings, $\epsilon$-admissibility interpolates between admissible covers with arbitrary ramifications over markings and relative stable maps.
As is explained in \cite{ACV},  the normalisation of a moduli space of admissible covers is a moduli space of stable twisted maps to $BS_n$, denoted by $\CK_{g,N}(BS_n,\mathsf m)$, where $\mathsf{m}$ is the number of simple ramifications which determines $\mathsf{h}$ by the Riemann--Hurwitz formula.  For the purposes of enumerative geometry (virtual intersection theory of moduli spaces),  the interpolation above can therefore be equally viewed in the following way, 
\[\xymatrix{
\CK_{g,N}(BS_n,\mathsf m) \ar@{<-->}[r]|-{\epsilon } &\Mbar_{\mathsf h}(C_{g,N},n).}
\]
In fact, this point of view is more suitable for the analogy with quasimaps. 
\subsection{Higher-dimensional case}We can upgrade the setting even further by adding a map $f_{X}\colon P \rightarrow X$ for some target variety $X$. This leads to the study of $\epsilon$-admissibility of  tuples
\[(P,C,\mathbf x, f_X \! \times \! f_C),\]
which can be represented as a correspondence
\[
\begin{tikzcd}[row sep=small, column sep = small]
P \arrow[r, "f_{X}"] \arrow{d}[swap]{f_{C}} & X  & \\
(C,\mathbf{x}) & & 
\end{tikzcd}
\]
In this case, $\epsilon$-admissibility also takes into account the degree of the components of $P$ with respect to the map $f_X$, see Definition \ref{epsilonadm}. If $X$ is a point, we recover the setting discussed previously.

Let $\beta=(\gamma, \mathsf m) \in H_2(X,\BZ)\oplus \BZ$ be an \textit{extended} degree. For $\epsilon \in (0,1]$, we then define \[Adm_{g,N}^{\epsilon}(X^{(n)},\beta)\] to be the moduli space of tuples
\[(P,C,\mathbf x, f_X\times f_C),\] such that $\mathsf{br}(f)=\mathsf m$,\footnote{By a version of Riemann--Hurwitz formula, Lemma \ref{RHformula}, the degree of the branching divisor $\br(f)=\mathsf m$ and the genus $\mathsf h$ determine each other.} $g(C)=g$, $|\mathbf{x}|=N$, and the map $f_X \times f_C$ is $\epsilon$-admissible of degree $(\gamma,n)$. The notation is slightly misleading, as $\epsilon$-admissible maps are not exactly maps to $X^{(n)}$. However, it is justified by the analogy with quasimaps and is more natural with respect to our notion of degree of $\epsilon$-admissible maps (see also Section \ref{Relation1}). 

As in the case of $X$ is a point, we obtain the following description of these moduli spaces for extremal values of $\epsilon$,
\begin{align*}
&\Mbar_{\mathsf h}(X\times C_{g,N}, (\gamma,n))= Adm_{g,N}^{0^+}(X^{(n)},\beta),\\
&\CK_{g,N}([X^{(n)}],\beta)\xrightarrow{\rho}  Adm_{g,N}^{1}(X^{(n)},\beta),
\end{align*}
such that the map $\rho$ is a virtual normalisation in the sense of the diagram (\ref{normalisation}), which makes two spaces equivalent from the perspective of enumerative geometry, see \cite{AGV} for the construction of $\CK_{g,N}([X^{(n)}],\beta)$ and \cite{BG} for the notion of extended degree of maps to orbifolds that we use in this work. We therefore get an interpolation,
\[\xymatrix{
\CK_{g,N}([X^{(n)}],\beta) \ar@{<-->}[r]|-{\epsilon } &\Mbar_{\mathsf h}(X\times C_{g,N}, (\gamma,n)),}
\]
which is completely analogous to (\ref{qmwall}).  
\subsection{Wall-crossing}  The invariants of $\Mbar_{\mathsf h}(X\times C_{g,N}, (\gamma,n))$ that can be related to orbifold invariants of $\CK_{g,N}([X^{(n)}],\beta)$ are the \textit{relative} GW invariants taken with respect to the markings of the target curve $C$. More precisely, for all $\epsilon$, there exist natural evaluation morphisms 
\begin{align*}
\ev_i\colon& \vec{A}dm_{g,N}^{\epsilon}(X^{(n)},\beta)\rightarrow \vec{\CI}X^{(n)}, \ i=1,\dots, N, 
\end{align*}
where $\vec{A}dm_{g,N}^{\epsilon}(X^{(n)},\beta)$ is defined by putting a standard order\footnote{This is an order which respects the ramification degrees of points.} on ramification points, and $\vec{\CI}X^{(n)}$ is an \textit{ordered} version of the inertia stack $\CI X^{(n)}$, given by connected components of $\CI X^{(n)}$ without the action of associated centraliser subgroups, see Section \ref{Section_Inertiastack}.   We define 
\begin{multline*}
\langle \psi_1^{m_{1}}\gamma_{1}, \dots, \psi_N^{m_{N}}\gamma_{N} \rangle^{\epsilon}_{g,N,\beta}:= \\
 \frac{1}{ \prod_i |\Aut(\mu^i)|}\int_{[ \vec{A}dm_{g,N}^{\epsilon}(X^{(n)},\beta)]^{\mathrm{vir}}}\prod^{i=N}_{i=1}\psi_{i}^{m_{i}} \ev^{*}_{i}(\gamma_{i}),
\end{multline*}
where $\gamma_{i}$ are classes in the cohomology $H^*(\vec{\CI}X^{(n)},\BQ)$, $\psi_i$ are $\psi$-classes associated to the markings of the target curves, and  $|\Aut(\mu^i)|$ are the factors introduced by the ordering of ramification points. By Section \ref{comp}, these invariants specialise to orbifold GW invariants associated to a moduli space $\CK_{g,N}([X^{(n)}],\beta)$ and relative GW invariants associated to a moduli space $\Mbar_{\mathsf h}(X\times C_{g,N}, (\gamma,n))$ for corresponding values of $\epsilon$. 

To relate invariants for different values of $\epsilon$, we use the master space technique developed by Zhou in \cite{YZ} for the purposes of the quasimap theory. We establish the properness of the master space in our setting in Section \ref{master}, following the strategy of Zhou. To state compactly the wall-crossing formula, we define 
\[F^{\epsilon}_{g}(\mathbf{t}(z))=\sum^{\infty}_{N=0}\sum_{\beta}\frac{q^{\beta}}{N!}\langle \mathbf{t}(\psi_1), \dots, \mathbf{t}(\psi_N) \rangle^{\epsilon}_{g,N,\beta},\]
where $\mathbf{t}(z) \in H^{*}(\vec{\CI}X^{(n)},\BQ)[z]$ is a generic element, and the unstable terms are set to be zero. There exists an element 
\[I_+(z) \in H^{*}(\vec{\CI}X^{(n)},\BQ)[z]\otimes_\BQ \BQ[\![q^{\beta}]\!],\] 
defined in Section \ref{graphspaceSym} as a truncation of the $I$-function. The $I$-function is in turn defined via the virtual localisation on the space of stable maps to $X\times \p^1$ relative to $X\times\{\infty\}$. The element $I_+(z)$ provides the change of variables, which relates generating series for extremal values of $\epsilon$. 
\begin{thmn}  For all $g\geq 1$, we have
\[F^{0^+}_{g}(\mathbf{t}(z))=F^{1}_{g}(\mathbf{t}(z)+I_+(-z)).\]
For $g=0$, the same equation holds modulo constant and linear terms in $\mathbf{t}$. 	
\end{thmn}
The change of variables above is the consequence of a wall-crossing formula  across each wall between extremal values of $\epsilon$, see Theorem \ref{wallcrossingSym}.
\subsection{Applications}

\subsubsection{The Square} For a del Pezzo surface $S$, we compute the wall-crossing invariants in Section \ref{del Pezzo}. A computation for analogous quasimap wall-crossing invariants is given in \cite[Proposition 8.6]{N}.

The wall-crossing invariants can easily be shown to satisfy \textsf{PT/GW}. Hence when both quasimap and \textsf{GW/H} wall-crossings are applied, \textsf{CRC} becomes equivalent to  \textsf{PT/GW}. For precise statements of both in this setting, we refer to Section \ref{qausiadm}. This is expressed in terms of the Square of theories in Figure \ref{square}. 

In \cite{PP},  \textsf{PT/GW} is established for $S\times \p^1$ relative to $S\times\{0,1,\infty\} \subset S\times \p^1$, if $S$ is toric. Alongside with  \cite{PP},  the square therefore gives us the following result. 
\begin{thmn}If $S$ is a toric del Pezzo surface, $g=0$ and $N=3$, then $\mathsf{CRC}$ (in the sense of \cite{BG}) holds for $S^{[n]}$ for all $n\geq 1$ and in all classes. 
\end{thmn}

Previously, the theorem above was established for $n = 2$ and $S = \p^2$ in \cite[Section 6]{W}; for an arbitrary $n$ and an arbitrary toric surface, but only for exceptional curve classes, in \cite{Che}; for an arbitrary $n$ and a simply connected $S$, but only for exceptional curve classes  and in the sense of \cite{YR}, in \cite{LQ}. We believe that with a little bit of effort, all of the results above can be given a more natural proof by reducing them to $\mathsf{PT/GW}$ for $S\times \p^1$ by means of our wall-crossings, as $\mathsf{PT/GW}$ is a more computationally accessible side.

 If $S = \BC^2$, \textsf{CRC} was proved for all genera and any number of markings on the level of cohomological field theories in \cite{PT19b}. If $S = \CA_n$, it was proved in genus-0 case and for any number of markings in \cite{CCIT} in the sense of \cite{CIR}.  For surfaces with $\mathrm{c}_1(S)=0$, $\mathsf{CRC}$ was established in \cite{FG} with representation-theoretic methods, as GW invariants vanish in this case.  \textsf{CRC}  was also proved for resolutions other than those that are  of the Hilbert--Chow type. The list is too long to mention them all.

The result is very appealing, because
the underlying cohomologies with classical multiplications are not isomorphic as rings for surfaces with $\mathrm{c}_1(S) \neq 0$,  but the associated genus-0 3-point invariants are equal up to a change of variables. In
particular, the classical multiplication on $H_{\mathrm{orb}}^*(S^{(n)},\BC)$ is a non-trivial quantum deformation of the classical multiplication on $H^*(S^{[n]},\BC)$. See Section \ref{qcoh} for more details.

We want to stress that $\textsf{CRC}$ is a more geometric side of the story than $\textsf{PT/GW}$, because it relates theories which are closer to each other. Moreover, as \cite{BG} points out and Section \ref{theremark} elaborates, the change of variables,
\begin{equation} \label{eiu}
	y=-e^{iq},
\end{equation}
arises naturally due to the following features of $\textsf{CRC}$,
\begin{itemize} 
	\item[(i)] analytic continuation of generating series from $0$ to $-1$,
	\item[(ii)] factor $i=\sqrt{-1}$ in the identification of cohomologies of $S^{[n]}$ and $S^{(d)}$,
	\item[(iii)] the divisor equation in $\mathsf{GW}(S^{[n]})$,
	\item[(iv)] failure of the divisor equation in $\mathsf{GW}_{\mathrm{orb}}([S^{(n)}])$.
\end{itemize} 
More precisely, (i) is responsible for the minus sign in (\ref{eiu}); (iii) and (iv) are responsible for the exponential; (ii) is responsible for $i$ in the exponential.  More conceptual view on \textsf{CRC} is presented in the works of Iritani, for example, \cite{I}. 

\subsubsection{LG/CY vs CRC} \label{LG/CY} We will now draw certain similarities between \textsf{CRC} and  Landau-Ginzburg/Calabi-Yau correspondence (\textsf{LG/CY}),  which are  illustrated in Table \ref{table}. For all details and notation on \textsf{LG/CY}, we refer to \cite{CIR}. For this discussion and the one in Section \ref{crepant}, it is more convenient to allow $\epsilon$ to be a negative number instead, $\epsilon \in [-1, 0) \subset \BR_{<0}$, in the definition of $\epsilon$-admissibility. This is to distinguish $\epsilon$-admissibility from $\epsilon$-stability of quasimaps discussed in \cite{N}, where $\epsilon$ also lives in $\BR_{>0}$.  See Section \ref{crepant} for more details.

\textsf{LG/CY} consists of two types of correspondences:  A-model and B-model correspondences. The B-model correspondence is the statement of equivalence of two categories - matrix factorisation categories and derived categories.  While the A-model correspondence is the statement of equality of generating series of certain curve-counting invariants after an analytic continuation and a change of variables. Moreover, there exists a whole family of enumerative theories depending on a stability parameter $\epsilon \in \BR$. For $\epsilon \in \BR_{>0}$ it gives the theory of GIT quasimaps, while for $\epsilon \in \BR_{< 0}$ it gives the FJRW (Fan--Jarvis--Ruan--Witten) theory. The GLSM (Gauged Linear Sigma Model) theory, defined mathematically in \cite{HTY}, allows to unify quasimaps and the FJRW theory.  The analytic continuation occurs, when one crosses the wall at $\epsilon=0$.

In the case of \textsf{CRC} we have a similar picture. The B-model correspondence is given by an equivalence of categories, $\mathrm{D}^b(S^{[n]})$ and $\mathrm{D}^b([S^{(n)}])$. The A-model correspondence is given by an analytic continuation of generating series and the subsequent application of a change of variables, as it is stated in Section \ref{crepant}. There also exists a family of enumerative theories depending on a parameter $\epsilon \in \BR$. For $\epsilon \in \BR_{> 0}$, it is given by  $\epsilon$-stable quasimaps to a moduli space of sheaves, while for $\epsilon \in \BR_{< 0}$ it is given by $\epsilon$-admissible maps.  It would be interesting to know  if a unifying theory exists in this case (like the GLSM theory in  \textsf{LG/CY}).

The above comparison is not a mere observation about structural similarities of two correspondences. In fact, both correspondences are instances of the same phenomenon. Namely, in both cases there should exist \textit{K\"ahler moduli spaces}, $\CM_{\mathsf{LG/CY}}$ and $\CM_{\mathsf{CRC}}$,  such that two geometries in question correspond to two different cusps of these moduli spaces (e.g.\ $S^{[n]}$ and $[S^{(n)}]$ correspond to two different cusps of $\CM_{\mathsf{CRC}}$).  B-models do not vary across these moduli spaces, hence the relevant categories are isomorphic. On the other hand,  A-models vary in the sense that there exist non-trivial global quantum $D$-modules, $\CD_{\mathsf{LG/CY}}$ and $\CD_{\mathsf{CRC}}$, which specialise to relevant enumerative invariants around cusps. For more details on this point of view, we refer to \cite{CIR} in the case of $\mathsf{LG/CY}$, and to \cite{I2} in the case of $\mathsf{CRC}$.

\vspace{.1in}
\begin{table*}[h!]
	\[
	\begin{tabular}{ |c|c|c| } 
		\hline
		& B-model & A-model \\ \hline 
		$\mathsf{LG/CY}$ & $\mathrm{D^b}(X_W)\cong \mathrm{MF}(W)$ & $\mathsf{GW}(X_W) \xleftarrow{\epsilon< 0}\mid_{0}\xrightarrow{\epsilon>0} \mathsf{FJRW}(\BC^{n},W)$ \\ 
		\hline 
		$\mathsf{CRC}$ & $\mathrm{D^b}(S^{[n]})\cong \mathrm{D^b}([S^{(n)}])$ & $\mathsf{GW}_{\mathrm{orb}}([S^{(n)}]) \xleftarrow{\epsilon<0}\mid_{0}\xrightarrow{\epsilon>0} \mathsf{GW}(S^{[n]})$ \\ 
		\hline
	\end{tabular}
	\]
	\caption{\textsf{LG/CY} vs \textsf{CRC}}
	\label{table}
\end{table*}
\subsection{Acknowledgments}
First and foremost I would like to thank Georg Oberdieck for the supervision of my PhD, and for pointing out that ideas of quasimaps can be applied to orbifold symmetric products. I am grateful to Daniel Huybrechts for reading some parts of the present work and  Maximilian Schimpf for providing the formula for Hodge integrals. 

A great intellectual debt is owed to Yang Zhou for his theory of entangled tails, without which the wall-crossings would not be possible. 
\subsection{Notation and conventions}
We work over the field of complex numbers $\BC$. Given a variety $X$, by $[X^{(n)}]$ we denote a stacky symmetric product $[X^n/S_n]$ and by $X^{(n)}$ its coarse quotient.  We denote a Hilbert scheme of $n$-points by $X^{[n]}$. For a partition $\mu=(\mu_1,\dots, \mu_k)$ of an integer $n$, let $\ell(\mu)=k$ denote the length of $\mu$ and $\mathrm{age}(\mu)=n-\ell(\mu)$. For a possibly disconnected curve $C$, we define $g(C):=1-\chi(\CO_C)$.

We set 
$e_{\BC^*}(\BC_{\mathrm{std}})=z,$
where $\BC_{\mathrm{std}}$ is the weight 1 representation of $\BC^*$ on a vector space $\BC$.  Let $N$ be a semigroup and $\beta \in N$ be its generic element. By $\BQ[\![ q^\beta ]\!]$ we will denote the (completed) semigroup algebra 
$\BQ[\![ N]\!]$. In our case, $N$ will be various semigroups of effective curve classes.

\section{Moduli spaces of $\epsilon$-admissible maps} \label{sectionadm}
\subsection{Branching divisor} Let $X$ be a smooth projective variety over the field of complex numbers $\BC$, $(C,\mathbf{x})$ be a marked nodal curve, and let $P$ be a possibly disconnected  nodal curve. 
\begin{defn} \label{twisted}
	For a map 
	\[f=f_{X}\times f_{C} \colon P \rightarrow X\times (C,\mathbf{x}),\]	the
	data $(P, C, \mathbf{x},f)$ is called a \textit{pre-admissible} map, if
	\begin{itemize}
		\item[$\bullet$] $f_{C}$ is admissible over nodes and marked points,
		\item[$\bullet$] $f$ is non-constant on each connected component of $P$.  
	\end{itemize}
\end{defn}
We will refer to $P$ and $C$ as \textit{source} and \textit{target} curves, respectively. We define the automorphism group of $f$ to be the  group of automorphisms on the source and  target curves that fix $f$, 
\[\Aut(f):=\{(\phi_1,\phi_2)\in \Aut(P)\times \Aut(C,\mathbf{x}) \mid f \circ \phi_1 = \phi_2 \circ f \}.\]

Let $\mathbf{x}'\subset P$ be the preimages of the markings $\mathbf{x}$.  We have a naturally defined map between logarithmic cotangent bundles of $(C,\mathbf{x})$ and $(P,\mathbf{x}')$, 
\begin{equation} \label{themap}
f_C^*\Omega_{(C,\mathbf{x})}^{\mathrm{log}} \rightarrow \Omega_{(P,\mathbf{x}')}^{\mathrm{log}},  
\end{equation}
we remind the reader that $\Omega_{(C,\mathbf{x})}^{\mathrm{log}}$ is naturally isomorphic to the twisted dualizing sheaf $\omega_C(\mathbf{x})$, the same holds for $P$, see \cite[Proposition 1.13]{Kato}, and \cite[Section 7]{AC} for the properties of log cotangent complexes in general. By the pre-admissibility, one can readily verify that the map (\ref{themap}) is an isomorphism around nodes and marked point by passing to local charts.  Let us view the map (\ref{themap}) as a complex, such that $\Omega_{(P,\mathbf{x}')}^{\mathrm{log}}$ is in the degree 0.  Consider now the following perfect complex,
\[Rf_{C*} [f_C^*\Omega_{(C,\mathbf{x})}^{\mathrm{log}} \rightarrow \Omega_{(P,\mathbf{x}')}^{\mathrm{log}}] \in \mathrm{D}^b(C),\]
which  is supported at finitely many points, which we call \textit{branching} points.  By the pre-admissibility, they arise either due to ramification points away from marked points or contracted components of the map $f_{C}$.  Following \cite{FP}, to the complex above, 
we can associate an effective Cartier divisor, 
\[\br(f) \in \mathrm{Div}(C),\]
by taking the support of the complex weighted by its Euler characteristics. This divisor will be referred to as a \textit{branching divisor}. 

Let us give a more explicit expression for the branching divisor. 
Let $P_{\circ} \subseteq P$ be the maximal subcurve of $P$ which is contracted by the map $f_{C}$.   Let $P_{\bullet} \subseteq P$ be the complement of $P_{\circ}$, i.e., the maximal subcurve which is not contracted by the map $f_{C}$. By $\widetilde P_{\bullet}$ we denote its normalisation at the nodes which are mapped into a regular locus of $C$. Note that the restriction of $f_{C}$ to $\widetilde P_{\bullet}$ is a ramified cover, the branching divisor of which is therefore given by points of ramifications.  

By $\widetilde P_{\circ, i}$ we denote the connected components of the normalisation $\widetilde P_{\circ}$ and by $p_i\in C$ their images in $C$. Finally, let $N \subset P$ be the locus of nodal points which are mapped into regular locus of $C$.  Following \cite[Lemma 10, 11]{FP}, the branching divisor $\br(f)$ can be expressed as follows. 
\begin{lemma} \label{br} With the notation from above, we have
	\[ \br(f)=\br(f_{|\widetilde P_{\bullet}}) + \sum_i(2g(\widetilde P_{\circ, i})-2)[p_i]+2f_*(N).\]
\end{lemma}
\textit{Proof.} This follows directly from  \cite[Lemma 10, 11]{FP}, since $\br(f)$ is supported on the regular locus of $C$. 
\qed
\\

We fix $L \in \Pic(X)$, an ample line bundle on $X$, such that for all effective curve classes $\gamma \in H_2(X,\BZ)$,
\[\deg(\gamma)= \mathrm{c}_1(L) \cdot \gamma > 1.\]
If  $\mathrm{c}_1(L) \cdot \gamma=1$, then in the definition of $\epsilon$-admissibility below, a simple ramification and a rational contracted component whose degree is equal to 1 with respect to the line bundle $L$ are indistinguishable. This should be avoided, because we want the maps to have simple ramifications but not contracted components. 

Let $(P,C,\mathbf{x}, f)$
be a pre-admissible map. For a point $p\in C$, we define  
\[f^*L_p:=f_X^*L_{|f_{C}^{-1}(p)},\] 
and set $\deg(f^*L_p)=0$, if $f_{C}^{-1}(p)$ is 0-dimensional. For a component $C'\subseteq C$, let 
\[f^*L_{|C'}:=f_X^*L_{|f_{C}^{-1}(C')}.\]
Recall that a \textit{rational tail} of a curve $C$ is a component isomorphic to $\p^1$ with one special point (a node or a marked point). A \textit{rational bridge} is a component isomorphic to $\p^1$ with two special points.
\begin{defn} \label{epsilonadm} Let $\epsilon \in (0,1]$.
	A  pre-admissible map $f$ is $\epsilon$-admissible, if 
	\begin{itemize}
		\item[(i)] for all points $p \in C$, 
		\[ \mathrm{mult}_p(\br(f))+\deg(f^*L_p)\leq 1/\epsilon,\]
		\item[(ii)] for all rational tails $T \subseteq (C,\mathbf{x})$, 
		\[\deg(\br(f)_{|T})+\deg(f^*L_{|T})>1/\epsilon,\]
		
		\item[(iii)] the group of automorphisms of $f$ is finite, 
		\[ |\Aut(f)|<\infty.\] 
		\end{itemize}
\end{defn}

A family of $\epsilon$-admissible maps over a base scheme $B$ is given by two families of curves $P$ and $(C,\mathbf{x})$ over $B$ and an admissible map
\[f=f_{X}\times f_{C} \colon P \rightarrow X\times (C,\mathbf{x}),\]
whose fibers over geometric points of $B$ are $\epsilon$-admissible. An isomorphism of families 
\[\Phi=(\phi_1, \phi_2)\colon (P,C,\mathbf{x}, f) \cong (P',C',\mathbf{x}', f')\]
is given by the data of isomorphisms of the source and target curves 
\[(\phi_1, \phi_2) \in \mathrm{Isom}_B(P,P') \times \mathrm{Isom}_B((C,\mathbf{x}),(C',\mathbf{x}') ),\]
such that   
\[f' \circ \phi_1 = \phi_2 \circ f.\]
The condition of $\epsilon$-admissibility is  open in the space of pre-admissible maps, as explained in the next lemma. 

\begin{lemma} \label{open} The condition of $\epsilon$-admissibility is an open condition in the space of pre-admissible maps.
\end{lemma}
\textit{Proof.} The conditions of $\epsilon$-admissibility are constructable.  Hence we can use the valuative criteria for openness. Given a discrete valuation ring $R$ with a fraction field $K$, we therefore need to show that if a pre-admissible map 
\[(P, C, \mathbf{x},f)\]
is $\epsilon$-admissible at a closed fiber $\Spec \BC$ of  $\Spec R$, then it is  $\epsilon$-admissible at the generic fiber. In fact, each condition of $\epsilon$-admissibility is open. The condition (iii) is standard. Let us therefore address the first two conditions.  

First, let 
\[T \subseteq (C ,\mathbf x) \]
be a family of subcurves of $(C ,\mathbf x)$ such that the generic fiber $T_{|\Spec K}$ is a rational tail that does not satisfy the condition (ii). Then the central fiber $T_{|\Spec \BC}$ of $T$ will be a tree of rational curves, whose rational tail does not not satisfy the condition (ii), because the degree of both $\br(f)$ and $f_X^*L$ can only decrease on the rational tail of $T_{|\Spec \BC}$. This shows that the condition (ii) is open.

   Similarly, the multiplicity of the branching divisor and the degree of a contracted curve with respect of $L$ at a given point can only increase at the central fiber, leading us to the same conclusion regarding the condition (i) of $\epsilon$-admissibility. 
\qed


\begin{defn} Given an element
	\[ \beta=(\gamma, \mathsf m) \in H_2(X,\BZ)\oplus \BZ,\]
	we say that an  $\epsilon$-admissible map $f$
	is of degree $\beta$ to $X^{(n)}$, of genus $g$ with $N$ markings, if 
	\begin{itemize}
		\item[$\bullet$] $f$ is of degree $(\gamma,n)$ and $\deg(\br(f))=\mathsf{m}$, 
		\item[$\bullet$] $g(C)=g$ and $|\mathbf{x}|=N$.
	\end{itemize}
\end{defn}

We define
\begin{align*}
	Adm_{g,N}^{\epsilon}(X^{(n)}, \beta) \colon (Sch/ \BC)^{\circ} &\rightarrow Grpd \\
	B&\mapsto \{\text{families of $\epsilon$-admissible maps over }B\}
\end{align*}
to be the moduli space of $\epsilon$-admissible maps  to $X^{(n)}$ of degree $\beta$ and genus $g$ with $N$ markings. 
Following \cite[Section 3.2]{FP}, one can construct the universal branching divisor
\begin{equation} \label{globalbrtw}
	\br:  Adm_{g,N}^{\epsilon}(X^{(n)}, \beta) \rightarrow \FM_{g,N,\mathsf{m}}.
\end{equation}
The space $\FM_{g,N,\mathsf m}$ is an algebraic stack which parametrises triples 
\[(C,\mathbf{x},D),\]
where $(C, \mathbf{x})$ is a genus-$g$ curve with $N$ markings; $D$ is an effective divisor of degree $\mathsf m$ disjoint from markings $\mathbf{x}$. An isomorphism of triples is an isomorphism of curves which preserve markings and divisors.

  Moduli spaces $ Adm_{g,N}^{\epsilon}(X^{(n)}, \beta)$  admit a disjoint-union decomposition 
\begin{equation}\label{ramificationprofiles}
	Adm_{g,N}^{\epsilon}(X^{(n)}, \beta)= \coprod_{\underline\mu} Adm_{g,N}^{\epsilon}(X^{(n)}, \beta,\underline{\mu}),
\end{equation}
where $\underline \mu=(\mu^1,\dots, \mu^N)$ is a $N$-tuple of partitions of $n$ which specify ramification profiles of $f_{C}$ over the markings $\mathbf{x}=(x_1, 
\dots, x_N)$. The Riemann--Hurwitz formula extends to the case of pre-admissible maps in the following form. 

\begin{lemma} \label{RHformula} If $f\colon P\rightarrow (C,\mathbf{x})$ is a degree $n$ pre-admissible map with ramification profiles $\underline{\mu}=(\mu^1,\dots, \mu^N)$ at the markings $\mathbf x \subset C$, then
	\[2g(P)-2=n\cdot (2g(C)-2)+\deg(\br(f))+ \sum_i\mathrm{age(\mu^i)},\]
	where $\mathrm{age(\mu)}=n-\ell(\mu)$.
	
\end{lemma} 
\textit{Proof.} Using Lemma \ref{br} and the standard Riemann--Hurwitz formula, one can readily check that the above formula holds for pre-admissible maps. \qed 

\subsection{Properness}
We now establish the properness of $ Adm_{g,N}^{\epsilon}(X^{(n)}, \beta)$, starting with the following result. 
\begin{prop} \label{DM1} A moduli space $ Adm_{g,N}^{\epsilon}(X^{(n)}, \beta)$ is a quasi-separated Deligne--Mumford stack of finite type. 
\end{prop}
\textit{Proof.} First, $ Adm_{g,N}^{\epsilon}(X^{(n)}, \beta)$ is an algebraic stack due to the representability of mapping stacks and Lemma \ref{open}.  Second, it is Deligne--Mumford, since automorphisms of maps are unramified by the condition (iii) of $\epsilon$-admissibility. It remains to show that it is of finite type and quasi-separated. 

   By the $\epsilon$-admissibility,
the map $\br$ factors through a quasi-separated substack of finite type. Indeed, $(C,\mathbf{x},\br(f))$ is not stable (i.e., has infinitely many automorphisms), if one of the following holds:
\begin{itemize} 
	\item[(i)] there is a rational tail $T \subseteq (C,\mathbf{x})$, such that $\mathrm{supp}(\br(f)_{|T})$ is at most a point,
	\item[(ii)] there is a rational bridge $B \subseteq (C,\mathbf{x})$, such that $\mathrm{supp}(\br(f)_{|B})$ is empty.
\end{itemize}
Up to a change of coordinates, the restriction of $f_{C}$ to $T$ or $B$ must of the form
\begin{equation} \label{constantmaps}
	z^{\underline n}\colon (\sqcup^k_{i=1} \p^1)\cup P' \rightarrow \p^1.
\end{equation} 
 Let us clarify the notation of (\ref{constantmaps}).  The curve $\sqcup^k_{i=1} \p^1$ is disjoint union of $k$ distinct $\p^1$. Over a rational tail $T$, a possibly disconnected marked nodal curve  $(P',\mathbf p)$ is attached via markings to the disjoint union  $\sqcup^k_{i=1} \p^1$ at points $0\in \p^1$  of connected components of the disjoint union; $P'$ is contracted to $0\in \p^1$  in the target $\p^1$; while on the $i$-th  $\p^1$ in the disjoint union, the map is given by $z^{n_i}$ for $\underline n=(n_1, \dots, n_k)$. Over a rational bridge $B$, no curve $P'$ is attached ($P'$ is an empty curve). 

The fact that the restriction of $f_{C}$ is given by a map of such form can be seen as follows. The conditions (i) or (ii) imply that the restriction of $f_{C}$ to $T$ or $B$ has at most two\footnote{Remember that branching might also be present at the nodes.} branching points, which in turn implies that the source curve must be $\p^1$ by the Riemann--Hurwitz formula. A map from $\p^1$ to itself with two ramifications points is given by $z^k\colon \p^1 \rightarrow \p^1$ up to a change of coordinate. For a rational tail $T$, there might also be a contracted component $P'$ attached to the ramification point.

In the case of (ii),  $\epsilon$-admissibility then requires
\[\deg(f^*L_{|B})>0.\]
While in the case of (i), 
\[\deg(\br(f)_{|T})=\mathrm{mult}_p(\br(f))\] for the unique point $p \in T$ which is not a node. Hence, in this case, $\epsilon$-admissibility requires 
\[\deg(f^*L_{|T})-\deg(f^*L_p)>0.\]
Since we fixed the class $\beta$, the conclusions above bound the number of components $T$ or $B$ by $\deg(\gamma)$. Hence the image of $\br$ is contained in a quasi-compact substack of $\FM_{g,N,\mathsf m}$, which is therefore quasi-separated and of finite type, because $\FM_{g,N,\mathsf m}$ is quasi-separated and locally of finite type. 

The branching-divisor map $\br$ is of finite type and quasi-separated, since a base change of $\br$ with respect to a base scheme $B \rightarrow \FM_{g,N,\mathsf m} $ is a subspace of the moduli space of admissible stable maps to $X\times C$ for a family of nodal curves $C \rightarrow B$. We therefore conclude that the moduli space $ Adm_{g,N}^{\epsilon}(X^{(n)}, \beta)$ is of finite type and quasi-separated itself, because $\br$ is of finite type, quasi-separated and factors through a quasi-separated substack of finite type. 
\qed
\begin{lemma} \label{contraction}
	Given a pre-admissible map $(P,C,\mathbf{x}, f)$.
	Let $(P',C',\mathbf{x}', f')$ be given by contraction of a rational tail $T\subseteq (C,\mathbf x)$ and stabilisation of the induced map 
	\[f\colon P \rightarrow X\times  C'.\] 
	Let $p \in C'$ be the image of contraction of $T$. Then the following holds
	\[\deg(\br(f)_{|T})+\deg(f^*L_{|T})=\mathrm{mult}_{p}(\br(f'))+\deg(f'^*L_{p}).\]
\end{lemma}
\textit{Proof.}
By Lemma \ref{RHformula}, 
\[2g(P_{|T})-2=-2n+\deg(\br(f))+n-\ell(p),\]
where $\ell(p)$ is the number of points in fiber above $p$,
from which it follows that
\begin{align*}
	\deg(\br(f))&=2g(P_{|T})-2+2n-n+\ell(p) \\
	&=2g(P_{|T})-2+n+\ell(p).
\end{align*}
By Lemma \ref{br}, 
\begin{align*}
	\mathrm{mult}_p(\br(f))&=2g(P_{|T})-2+2\ell(p)+n-\ell(p)\\
	&=2g(P_{|T})-2+n+\ell(p).
\end{align*} 
It is also clear by definition, that 
\[\deg(f^*L_{|T})=\deg(f^*L_p),\]
the claim then follows. 
\qed 

\begin{defn} \label{modifiation}
	Let $R$ be a discrete valuation ring. Given a pre-admissible map 
	$(P,C,\mathbf{x}, f)$ over $\Spec R$. A \textit{modification} of $(P,C,\mathbf{x}, f)$ is a pre-admissible map
	$(\widetilde{P},\widetilde{C},\widetilde{\mathbf{x}},\widetilde f)$ over $\Spec R'$, such that 
	\[(\widetilde{P},\widetilde{C},\widetilde{\mathbf{x}}, \widetilde{f})_{|\Spec K'} \cong (P,C,\mathbf{x}, f)_{|\Spec K'},\]
	where $R'$ is a finite extension of $R$ with a fraction field $K'$. 
	
\end{defn}

A modification of a family of curves $C$ over a discrete valuation ring is given by three operations:
\begin{itemize}
	\item blow-ups of the central fiber of $C$,
	\item contractions of rational tails and rational bridges of the central fiber of $C$,
	\item base changes with respect to finite extensions of discrete valuation rings.
\end{itemize}
A modification of a pre-admissible map is therefore given by an appropriate choice of three operations above applied to both target and source curves, such that the map $f$ can be extended as well. 
\begin{thm} \label{properness} 
	The moduli spaces $ Adm_{g,N}^{\epsilon}(X^{(n)}, \beta)$ are proper Deligne-Mumford stacks. 
\end{thm}

\textit{Proof.} We will now use the valuative criteria of properness for quasi-separated Deligne--Mumford stacks. Let 
\[(P^{*},C^{*},\mathbf{x}^*, f^*) \in  Adm_{g,N}^{\epsilon}(X^{(n)}, \beta)(K)\]
be a family of $\epsilon$-admissible maps over the fraction field $K$ of a discrete valuation ring $R$. The strategy of the proof is to separate $P^*$ into two components $P^*_{\circ}$ and $P^*_{\bullet}$, the components contracted and not contracted by $f^*_{C^*}$, respectively (as it was done for Lemma \ref{br}). We then take a limit of $f^*_{|P^*_{\bullet}}$ preserving it as a cover over the target curve, and a limit of $f^*_{\mid P^*_\circ}$ as a stable map to $X$. Next, we glue the two limits back and perform a series of modifications to get rid of points or rational tails that do not satisfy $\epsilon$-admissibility.   
\\

\textit{Existence, Step 1.}  Let 
\[(P^{*}_{\circ}, \mathbf{q}_\circ^*) \subseteq P^{^*} \]
be the maximal subcurve contracted by $f^{^*}_{C^*}$, the markings $\mathbf{q}^{*}_{\circ}$ are given by the nodes of $P^{*}$ disconnecting $P^{*}_{\circ}$ from the rest of the curve. By
\[(P^{*}_{\bullet}, \mathbf{q}_{\bullet}^*) \subseteq P^{^*} \]
we denote the complement of $P^{*}_{\circ}$ with similar markings. 
Let 
\[(\widetilde{P}^{*}_{\bullet}, \mathbf{t}_1^*, \mathbf{t}_2^*)\]
be the normalisation of $P^{*}_{\bullet}$ at nodes which are mapped by $f^{*}_{C^*}$ to the regular locus of $C^{*}$, the markings $\mathbf{t}_1^*$ and $\mathbf{t}_2^*$ are given by the preimages of those nodes. The induced map 
\[\tilde{f}^{*}_{\bullet,C^*}\colon \widetilde{P}^{*}_{\bullet} \rightarrow (C^*,\mathbf{x}^*,\tilde{\mathbf{x}}^*)\] 
is an admissible cover, where $\tilde{\mathbf{x}}^*$ is an extra set of markings that accounts for additional ramifications of the map $\tilde{f}^{*}_{\bullet,C^*}$ allowed by the $\epsilon$-admissibility.  Consider now the  map \[
\tilde{f}^{*}_{\bullet}\colon \widetilde{P}^{*}_{\bullet} \rightarrow X\times (C^*,\mathbf{x}^*,\tilde{\mathbf{x}}^*),
\] it can be viewed as an element of  $\CK_{g,N+\tilde N}([X^{(n)}])(K)$ (see Section \ref{Relation1} for more details on how to do this), where $\tilde N=|\tilde{\mathbf{x}}^*|$.  Using the properness of $\CK_{g,N+\tilde N}([X^{(n)}])$ from \cite{AGV}, we  can extend this map to $\Spec R$, possibly after a finite base change and taking the coarse curves after the extension, 
\[
\tilde{f}_{\bullet}\colon \widetilde{P}_{\bullet} \rightarrow X\times (C,\mathbf{x},\tilde{\mathbf{x}}).
\] 
By the construction of $\CK_{g,N+\tilde N}([X^{(n)}])$, this extension has the property that the map $\tilde f_{\bullet | C}$ is an admissible cover ramified over $\mathbf{x}$, $\tilde{\mathbf{x}}$, such that the ramification profiles of generic and central fibers agree.\footnote{Such extension can also be constructed differently by firstly taking the extension of the admissible cover over $C^*$ and then extending the map to $X$ by blowing up $C$. }
We then reintroduce back the markings $(\mathbf{q}_\bullet^*, \mathbf{t}_1^*, \mathbf{t}_2^*)$. Note that after extending them to $\widetilde{P}_{\bullet}$, they might intersect between themselves as well as with the ramification points above $\mathbf{x}$ or with the singular locus. If this is the case, we blow up both the source curve  $\widetilde{P}_{\bullet}$ and the target curve $C$ at the intersecting markings and their images to separate them. This is always possible, because markings $(\mathbf{q}_\bullet^*, \mathbf{t}_1^*, \mathbf{t}_2^*)$ and ramification points above $\mathbf{x}^*$ are pairwise disjoint and contained in the smooth locus by construction.  Blow-ups will introduce chains of $\p^1$ in the central fibers of the source and the target. By the universality of blow-ups, the map $f_C$ extends,  such that its restriction over these chains is given by $z^k \colon \p^1 \rightarrow \p^1$  on each component, where $k$ will be determined by the ramification of the map at the markings before the blow-up;\footnote{This can be seen by passing to the formal neighbourhood of a point. The ramification points at the generic and central fibers agree, hence in a formal neighbourhood of a smooth point of the central fiber, the map between two families of curves therefore looks like $\BC[\![ x,y]\!]\rightarrow \BC[\![ x,y]\!]$, $x \mapsto x^k, y\mapsto y$. } while the chains are mapped trivially to $X$.  In particular, the pre-admissibility of maps is preserved. We take the minimal blow-up to avoid any unnecessary components (i.e., those components over which the map might have infinitely many automorphisms). Overall, after relabelling the resulting blown up curves,   we obtain a map over $\Spec R$,
\[
\tilde{f}_{\bullet}\colon (\widetilde{P}_{\bullet}, \mathbf{q}_\bullet, \mathbf{t}_1, \mathbf{t}_2) \rightarrow X\times (C,\mathbf{x},\tilde{\mathbf{x}}),
\] 
such that over $C$ it is an admissible cover ramified over $(\mathbf{x},\tilde{\mathbf{x}})$. Moreover, the degrees of ramifications are preserved over $(\mathbf{x},\tilde{\mathbf{x}})$,  and $(\mathbf{q}_\bullet, \mathbf{t}_1, \mathbf{t}_2)$ together with ramification points above $\mathbf{x}$ are pairwise disjoint. The group of automorphisms of this map is finite, because it was so in $\CK_{g,N+\tilde N}([X^{(n)}])$ and, by construction, all rational components introduced by blow-ups contain enough markings to stabilize them (otherwise, those markings were either disjoint at the central fiber or equal at the generic fiber which contradicts the assumptions for blow-ups). 

Now let 
\[f_{\circ,X} \colon  (P_{\circ},\mathbf{q}_{\circ}) \rightarrow X\] 
be the extension of $f_{\circ,X} \colon  (P^*_{\circ},\mathbf{q}^*_{\circ}) \rightarrow X$
over $\Spec R$ as a stable map. It exists, possibly after a finite base change, by the properness of the moduli space of stable marked maps to a projective target $X$.  Finally, we glue back $P_{\circ}$ and $P_{\bullet}$ at the markings $(\mathbf{q}_{\circ},\mathbf{q}_{\bullet})$, and reintroduce the nodes by gluing the markings $(\mathbf{t}_1,\mathbf{t}_2)$ to obtain a pre-admissible map 
\[f\colon P \rightarrow X\times (C, \mathbf{x}),\]
note that by construction the markings at which we glue the curves have the same image with respect to maps to $X$. Hence the gluing of curves together with maps to $X$ is possible, while $P_\circ$ is contracted by $f_C$. 
 Next, we perform a series of modifications to the map above to obtain an $\epsilon$-admissible map. 
\\

\textit{Existence, Step 2.} Let us analyse $(P,C,\mathbf{x}, f)$ in relation to the conditions of $\epsilon$-admissibility. 
\\

(i)  Assume there is a point in the central fiber $C_{| \Spec \BC}$ that does not satisfy the condition (i) of $\epsilon$-admissibility. By the $\epsilon$-admissibility of the generic fiber and the construction of   $(P,C,\mathbf{x}, f)$, this is possible only if the central fibers of markings $(\mathbf{q}_\bullet, \mathbf{t}_1, \mathbf{t}_2)$ are contained in the same fiber of $\tilde{f}_{\bullet}$ in the way that $\epsilon$-admissibility is violated for $f$, i.e., the contribution of contracted components and nodes exceed the one of the generic fiber for which $\epsilon$-admissibility is satisfied.\footnote{For example, two contracted components which do not violate $\epsilon$-admissibility at the generic fiber might end up in the same fiber of the map at the central fiber of the family violating $\epsilon$-admissibility. In this case, we blow up the target in the image point of those contracted components, and the source at the points of attachment of the contracted components. The result is that the source have additional two rational tails attached to it with contracted components, while the target has one rational tail, such that the images of contracted components are distinct.} Indeed, since the ramifications of the central fiber and the generic fibers above $
\tilde{\mathbf{x}}$ are equal, they satisfy $\epsilon$-admissibility by the $\epsilon$-admissibility of the generic fiber, therefore the only excess contributions arise due to contracted components and nodes. 

The degrees of contracted components with respect to the line bundle $L$ at $(\mathbf{q}_\bullet, \mathbf{t}_1, \mathbf{t}_2)$ and their genera are left unchanged. Hence  whenever there is an excess contribution from the contracted components and nodes at the central fiber due to some markings from the set $(\mathbf{q}_\bullet, \mathbf{t}_1, \mathbf{t}_2)$, we can separate those markings  by further blowing up the markings themselves and their images in $C$. This is always possible, since by the $\epsilon$-admissibility of the generic fiber, the generic fibers of the markings must be contained in different fibers of the map. This will again introduce chains of $\p^1$ over which the map is given by $z^k \colon \p^1 \rightarrow \p^1$ for some $k$ determined by the ramification degree before the blow-up. The result of such blow-up is the reduction of the quantity 

\[\mathrm{mult}_p(\br(f))+\deg(f^*L_p).\]

By iterating this process, we obtain a pre-admissible map, such that the contribution of the contracted components  together with the ramifications over $\tilde{\mathbf{x}}$ and nodes does not exceed $1/\epsilon$ over all points in $C$. Also, no new rational tails which violate the condition (ii) of $\epsilon$-admissibility are created by Lemma \ref{contraction}.   We take a minimal blow-up with such property to again avoid unnecessary  rational components. 

This modification is completely analogues to the one needed for proving that Hassett's spaces \cite{Ha03} are proper. In fact, by assigning weights to points based on contracted components and ramifications, the only difference with \cite{Ha03} is that we have to modify both target and source curves to keep the maps admissible over nodes. 
 \\
 
(ii) If a rational tail $T \subseteq (C_{|\Spec \BC},\mathbf x_{|\Spec \BC})$ does not satisfy the condition (ii) of $\epsilon$-admissibility, we contract it.
This is possible because by the $\epsilon$-admissibility of the generic fiber, $T$ cannot be a limit of a rational tail in the generic fiber (i.e., it must be a negative self-intersection curve in the total space of the family).  By Lemma \ref{contraction},
the central fiber satisfies the condition (i) of $\epsilon$-admissibility at the image point of the contraction. We iterate this process until we get rid of all rational tails that do not satisfy the condition (ii) of $\epsilon$-admissibility. In the process we do not create any new points that do not satisfy the condition (i) of  $\epsilon$-admissibility. 
\\

(iii) Since the map we started with had a finite group of automorphisms, the result of the modifications described above also must have a finite group of automorphisms. Indeed, in the first modification, we take minimal blow-ups which do not introduce unnecessary rational components, while by construction, those rational components that were introduced must have enough ramifications and contracted components attached to them to be stable, as they were created in the process of separating the associated markings.  In the second modification, we contract rational components, this might create new rational tails or rational bridges. However, these new rational components are  stabilised  by either ramification points or contracted components attached to them (or over them, if we are talking about the target curves)  because they were stable before the contraction.   Overall, we constructed a map over $\Spec R$,
\[f\colon P \rightarrow X\times (C, \mathbf{x}),\]
which is $\epsilon$-admissible, showing the existence part. 
\\

\textit{Uniqueness.} Assume we are given two families of $\epsilon$-admissible maps over $\Spec R$,
\[(P_1,C_1,\mathbf{x}_1, f_1)\ \text{and} \ (P_2,C_2,\mathbf{x}_2, f_2),\] 
which are isomorphic over $\Spec K$. Possibly after a finite base change, there exists a family of pre-admissible maps 
\[(\widetilde{P},\widetilde{C}, \tilde{\mathbf{x}},\tilde{f})\] 
which dominates both families in the sense that there exists a commutative square
\begin{equation}\label{dominatin}
	\begin{tikzcd}[row sep=small, column sep = small] 
		\widetilde{P} \arrow[r, "\tilde{f}"] \arrow[d]&  X\times \widetilde{C} \arrow[d] & \\
		P_i \arrow[r,"f_i"]&  X\times C_i  &
	\end{tikzcd}
\end{equation}
It is constructed by taking a curve $\tilde{C}$ which dominates $C_1$ and $C_2$, and which is therefore given by blow-ups at the central fibers of $C_1$ or $C_2$. By blowing up the central fiber of source curves $P_1$ and $P_2$, we may extend the maps $f_1$ and $f_2$ to $\tilde{C}$. By the properness of moduli spaces of stable maps to a fixed target, these extensions must be unique and therefore isomorphic, since they are isomorphic at the generic fiber. We call the resulting family $(\widetilde{P},\widetilde{C}, \tilde{\mathbf{x}},\tilde{f})$.  We take a minimal family  with such property. 

The maps $\tilde{C}_i \rightarrow \tilde{C}$ in (\ref{dominatin}) are given by contraction of rational tails and rational bridges. In the case of rational bridges, this is possible only if these rational bridges are unstable, otherwise  $f_1$ and $f_2$ cannot be admissible over nodes (i.e., there must be a contracted component over nodes or the ramifications of the touching components are not equal). Since both families are $\epsilon$-admissible, if a rational bridge is contracted in the case of one family, it must be contracted in the case of another family too. 

Consider now a rational tail $T \subseteq (\tilde{C}_{| \Spec \BC}, \tilde{\mathbf{x}}_{|\Spec \BC})$,  contracted by a vertical map in (\ref{dominatin}). Then by Lemma \ref{contraction}, we have the equality
\[\deg(\br(f)_{|T})+\deg(L_{|T})=\mathrm{mult}_{p}(\br(f))+\deg(L_{p}),\]
where $p$ is the image of $T$ under the contraction in any of the two families. The equality above implies that $T$ cannot satisfy the condition (ii) of $\epsilon$-admissibility, if the image point $p$ satisfies the condition (i) of $\epsilon$-admissibility. Both families  are $\epsilon$-admissible by assumption, hence if a rational tail $T$ is contracted in the case of one family, then it must also be contracted in the case of another. By induction, this applies to trees of rational components contracted by the rightmost vertical maps in (\ref{dominatin}). Overall,  this implies that 
\[ 
C_1\cong \tilde{C} \cong C_2. \]
 By the separatedness of moduli spaces of stable maps to a fixed target,  it must be that 
\[ (P_1,C_1,\mathbf{x}_1, f_1) \cong (\widetilde{P},\widetilde{C}, \tilde{\mathbf{x}},\tilde{f}) \cong (P_2,C_2,\mathbf{x}_2, f_2),\] 
this shows the uniqueness part and finishes the proof. 
\qed
\subsection{Obstruction theory} \label{SectionObs}
The obstruction theory of $ Adm_{g,N}^{\epsilon}(X^{(n)}, \beta)$ is defined via the obstruction theory of relative maps in the spirit of \cite[Section 2.8]{GV} with the difference that we have a relative target geometry $X\times \FC_{g,N}/ \FM_{g,N}$.

 Consider the moduli space of admissible maps to the universal curve $\FC_{g,N}/ \FM_{g,N}$, 
\[ \FM_{\mathsf{h}}(\FC_{g,N}/ \FM_{g,N},n)\]
we endow with the obstruction theory from \cite[Section 2.8]{GV}. Then, by Lemma \ref{open}, there exists an obstruction theory on $Adm_{g,N}^{\epsilon}(X^{(n)}, \beta)$ relative to the stack  $\FM_{\mathsf{h}}(\FC_{g,N}/ \FM_{g,N},n)$ given by the complex
\[R\pi_*(f_X^*T_X),\]
where $f_X \colon \CP\rightarrow X$ is the  universal map to $X$ from the universl source curve, and $\pi \colon  \CP \rightarrow Adm_{g,N}^{\epsilon}(X^{(n)}, \beta)$ is the natural projection. Overall, we have a relative perfect obstruction theory

\[ \phi: R\pi_*(f_X^*T_X)^\vee \rightarrow \BL_{	 Adm_{g,N}^{\epsilon}(X^{(n)}, \beta)/\FM_{\mathsf{h}}(\FC_{g,N}/ \FM_{g,N},n)}.\]

The virtual fundamental class is defined via the virtual pullback of the virtual fundamental class of $\FM_{\mathsf{h}}(\FC_{g,N}/ \FM_{g,N},n)$, see \cite{Ma}. The obstruction theory can also be defined via the moduli spaces of maps between twisted curve \cite{ACV} or maps between logarithmic curves \cite{KiL} instead of  $\FM_{\mathsf{h}}(\FC_{g,N}/ \FM_{g,N},n)$, both of which provide certain virtual normalisations of spaces of admissible maps. The associated virtual fundamental classes are equal (after applying pushforwards) by results in the aforementioned references. 

\subsection{Relation to other moduli spaces} \label{Relation1} Let us now relate the moduli spaces of $\epsilon$-admissible maps for the extremal values of $\epsilon \in (0,1]$ to more familiar moduli spaces.
\subsubsection{The case of $\epsilon=1$}\label{comparing} In this case, the first two conditions of Definition \ref{epsilonadm} are:
\begin{itemize}
	\item[(i)] for all points $p \in C$, 
	\[ \mathrm{mult}_p(\br(f))+\deg(f^*L_p)\leq 1,\]
	\item[(ii)] for all rational tails $T \subseteq (C,\mathbf{x})$, 
	\[\deg(\br(f)_{|T})+\deg(f^*L_{|T})>1.\]
\end{itemize}
Since the multiplicity and the degree take only integer values, by Lemma \ref{br} and the choice of $L$, there is only one possibility for which the condition (i) is satisfied. Namely,
$f_{C}$ does not contract any irreducible components and has only simple ramifications. 

To unpack the condition (ii), recall that a non-constant ramified map from a smooth curve to $\p^1$ has at least two ramification points; it has precisely two ramification points, if it is given by 
\begin{equation} \label{simple}
	z^{2}\colon \p^1 \rightarrow \p^1
\end{equation}
up to a change of coordinates. Hence \[\mathrm{mult}_p(\br(f))+\deg(f^*L_p)=1,\]
if and only if $f_C=z^2$ and $f_X$ is constant. In this case,  $|\mathrm{Aut}(f)|=\infty$. In light of the condition (iii) of $\epsilon$-admissibility, the condition (ii) is therefore automatically satisfied. 

We obtain that the data of a $\epsilon$-admissible map 
$(P,C,\mathbf{x}, f)$ for $\epsilon=1$ can be represented as follows, 
\[
\begin{tikzcd}[row sep=small, column sep = small]
	P \arrow[r, "f_{X}"] \arrow{d}[swap]{f_{C}} & X  & \\
	(C,\mathbf{x}, \mathbf{p}) & & 
\end{tikzcd}
\]
where $f_{C}$ is a degree $n$ admissible cover with arbitrary ramifications over the marking $\mathbf{x}$ and with simple ramifications over the unordered marking $\mathbf{p}=\br(f)$, such that $|\mathrm{Aut}(f)| <\infty$.
Hence  the moduli space $Adm_{g,N}^{1}(X^{(n)},  \beta)$ admits a projection from the moduli space of twisted stable maps with \textit{extended degree} (see \cite[Section 2.1]{BG} for the definition) to the orbifold $[X^{(n)}]$, constructed in \cite{AGV,ACV},
\begin{equation} \label{projection}
	\rho \colon \CK_{g,N}([X^{(n)}], \beta) \rightarrow 	 Adm_{g,N}^{1}(X^{(n)}, \beta),
\end{equation}
which is given by passing from twisted curves to their coarse  spaces.	Indeed, there is a correspondence between $S_n$-torsors and degree $n$ \'etale covers of twisted curves given by associating to a $S_n$-torsor $\CP'$ the quotient $[\CP'/S_{n-1}]$, and, conversely, to a cover $\CP$ the $n$-fold product $\CP\times_\CC \CP \dots \times_\CC \CP$ without the big diagonal, i.e., without the locus where at least two points coincide. Equivalently, this follows from the fact that the group of local automorphisms of both types of objects is $S_n$.  Hence an  element of $ \CK_{g,N}([X^{(n)}], \beta)$ can be represented as follows,
\[
\begin{tikzcd}[row sep=small, column sep = small]
	\CP \arrow[r, "f_{X}"] \arrow{d}[swap]{f_{\CC}} & X  & \\
	(\CC,\mathbf{x}, \mathbf{p}) & & 
\end{tikzcd}
\]
where $f_{\CC}$ is a representable degree $n$ \'etale cover over twisted marked curve $(\CC,\mathbf{x}, \mathbf{p})$. The additional set of markings $\mathbf{p}$ is unordered, over these markings the map $f_{\CC}$ must have simple ramifications after passing to coarse moduli spaces. The map $f_{X}$ has to be fixed by only finitely many automorphisms of the cover $f_{\CC}$. Passing to coarse moduli space, the above data becomes the data of a $\epsilon$-admissible map for $\epsilon=1$. The morphism $\rho$ therefore associates a coarse curve to a stacky one. The virtual fundamental classes of two  moduli spaces are related by the next result. 
\begin{lemma} \label{fc} We have
	\begin{equation*} 
		\rho_*[\CK_{g,N}([X^{(n)}],\beta)]^{\mathrm{vir}} 
		= [ Adm_{g,N}^{1}(X^{(n)}, \beta)]^{\mathrm{vir}}.
	\end{equation*}
\end{lemma}
\textit{Proof.} Let $\mathfrak{K}_{g,N}(BS_n, \mathsf m)$ be the moduli stacks of twisted maps to $BS_n$ (not necessarily stable) and $\mathfrak{Adm}_{g, \mathsf h, n, N}$ be the moduli stack of admissible covers (again not necessarily stable). By the discussion above and the fact that a map to scheme from a stacky curve factors through the coarse curve, there exists the following pull-back diagram,
\begin{equation} \label{normalisation}
	\begin{tikzcd}[row sep=small, column sep = small]
		\CK_{g,N}([X^{(n)}],\beta) \arrow[d,"\pi_1"] \arrow[r,"\rho"] &  Adm_{g,N}^{1}(X^{(n)}, \beta) \arrow[d,"\pi_2"] & \\
		\mathfrak{K}_{g,N}(BS_n, \mathsf m) \arrow[r] & \mathfrak{Adm}_{g, \mathsf h, n, N}
	\end{tikzcd}
\end{equation}
such that the relative obstruction theories are compatible, as they are both given by the complex  $R\pi_*(f_X^*T_X)$.
The bottom arrow is a normalisation map by  \cite[Proposition 4.2.2]{ACV}, therefore it is of degree $1$.  By \cite[Theorem 5.0.1]{Co}, we therefore obtain the claim for virtual fundamental classes given by the relative obstruction theories, 
\begin{multline} \label{compatability}
	\rho_*[\CK_{g,N}([X^{(n)}],\beta)/\mathfrak{K}_{g,N}(BS_n, \mathsf m)]^{\mathrm{vir}} \\= [ Adm_{g,N}^{1}(X^{(n)}, \beta)/\mathfrak{Adm}_{g, \mathsf h, n, N}]^{\mathrm{vir}}.
\end{multline}
The moduli space $\mathfrak{K}_{g,N}(BS_n, \mathsf m)$ is smooth and connected, and it is the normalisation of $\mathfrak{Adm}_{g, \mathsf h, n, N}$, this implies that associated virtual fundamental classes of both spaces are equal to fundamental classes. 
Using virtual pull-backs of \cite{Ma}, one can therefore express the virtual fundamental classes given by absolute perfect obstruction theories as follows,
\begin{align*}[ Adm_{g,N}^{1}(X^{(n)}, \beta)]^{\mathrm{vir}}&=\pi_2^{!}[\mathfrak{Adm}_{g, \mathsf h, n, N}] \\
	&=[ Adm_{g,N}^{1}(X^{(n)}, \beta)/\mathfrak{Adm}_{g, \mathsf h, n, N}]^{\mathrm{vir}},
\end{align*}
the same applies to $\CK_{g,N}([X^{(n)}],\beta)$. Hence we obtain that 
\begin{equation*} 
	\rho_*[\CK_{g,N}([X^{(n)}],\beta)]^{\mathrm{vir}} 
	= [ Adm_{g,N}^{1}(X^{(n)}, \beta)]^{\mathrm{vir}},
\end{equation*}
this finishes the proof. 
\qed

\subsubsection{The case of $\epsilon=0^+$} By the first two conditions of Definition \ref{epsilonadm}, the map $f_{C}$ can have arbitrary ramifications and contracted components of arbitrary genera (more precisely, the two are only restricted by $n$, $g$, $N$ and $\beta$). In conjunction with other conditions of Definition \ref{epsilonadm} we therefore obtain the following identification of moduli spaces
\begin{equation} \label{compatability2}
	Adm_{g,N}^{0^+}(X^{(n)}, \beta) = \Mbar_{\mathsf{h}}(X\times C_{g,N},(\gamma,n)),
\end{equation}
where the space on the right is the moduli space of relative stable maps with possibly disconnected sources to the relative geometry \[X\times C_{g,N} \rightarrow \Mbar_{g,N},\] 
where $C_{g,N} \rightarrow \Mbar_{g,N}$ is the universal curve and where the markings play the role of relative divisors; curves are allowed to sprout rational bridges for the sake of properness and existence of perfect obstruction theory. Instead of fixing the genus of source curves, we fix the degree $\mathsf{m}$ of the branching divisor.  At each component $ Adm_{g,N}^{0^+}(X^{(n)}, \beta,\underline{\mu})$ of the decomposition (\ref{ramificationprofiles}), the genus of the source curve and the degree of the branching divisor are related by Lemma \ref{RHformula}.

The obstruction theories of two moduli spaces are equal, since the obstruction theory of the space $ Adm_{g,N}^{0^+}(X^{(n)}, \beta)$ was defined via the obstruction theory of relative stable maps. 
\subsection{Inertia stack} \label{Section_Inertiastack}
The inertia stack of a symmetric product can be defined as follows, 
\[\CI X^{(n)}=\coprod_{[g]}[X^{n,g}/C(g)],\] 
where the disjoint union is taken over conjugacy classes $[g]$ of elements of $S_n$, $X^{n,g}$ is the fixed locus of $g$ acting on $X^n$ and $C(g)$ is the centraliser subgroup of $g$. Recall that conjugacy classes of elements of $S_n$ are in one-to-one correspondence with partitions $\mu$ of $n$. Let us express a partition $\mu$ in terms of repeating parts and their multiplicities,
\[\mu=(\underbrace{ 1, \hdots, 1}_{m_1}, \hdots, \underbrace{s, \hdots, s}_{m_s}).\]
We define 
\begin{equation} \label{centraliser}
	C(\mu):= \prod^{s}_{t=1} C_{\eta_t} \wr S_{m_t},
\end{equation}
here $C_{\eta_t}$ is a cyclic group and $``\wr"$ is a \textit{wreath product} defined as 
\[C_{\eta_t} \wr S_{m_t}:=C_{\eta_t}^{\Omega_t}\rtimes S_{m_t},\]
where $\Omega_t=\{1,2,\dots,m_t\}$; $S_{m_t}$ acts on $C_{\eta_t}^{\Omega_t}$ by permuting the factors. There exist two natural subgroups of $C(\mu)$ 
\begin{equation} \label{autmu}
	\Aut(\mu):=\prod^{s}_{t=1} S_{m_t} \ \text{ and } \ N(\mu):=\prod^{s}_{t=1} C_{\eta_t}^{\Omega_t}
\end{equation}
as the notation suggests, $\Aut(\mu)$ coincides with the automorphism group of the partition $\mu$. On the other hand, $|N(\mu)|=\prod_j \mu_j$. By construction, we have the following sequence, 
\begin{equation} \label{split}
	1 \rightarrow N(\mu) \rightarrow C(\mu) \rightarrow \Aut(\mu) \rightarrow 1.
\end{equation}
In particular, we obtain that 
\[ \mathfrak{z}(\mu):=|\Aut(\mu)|\cdot \prod_j \mu_j=|C(\mu)|. \]

We put the \textit{standard order} on the partition, 
\[
\mu_i\geq \mu_j \iff j\geq i.
\]
  Viewing a partition $\mu$ as a partially ordered set, we define $X^{\mu}$ as the self-product of $X$ over the set $\mu$. In particular, 
\[X^{\mu}\cong X^{\ell (\mu)},\] 
where $\ell(\mu)$ is the length of the partition $\mu$. The group $C(\mu)$ acts on $X^\mu$ as follows. The products of cyclic groups $C_{\eta_t}^{\Omega}$ acts trivially on corresponding factors of $X^\mu$, while $S_{m_t}$ permutes the factors corresponding to the same part $\eta_t$. These actions are compatible with the wreath product. 

Given an element $g\in S_n$ in a conjugacy class corresponding to a partition $\mu$, we have the following identifications
\[C(g) \cong C(\mu)\ \text{and} \ X^{n,g}\cong X^\mu,\] 
such that the  group actions match.  
With the notation introduced above, the inertia stack can be re-expressed as follows, 
\begin{equation} \label{Inertia}
	\CI X^{(n)}=\coprod_{\mu}[X^{\mu}/C(\mu)]. 
\end{equation}
Motivated by the considerations above, we define the ordered inertia stack, 
\[ \vec{\CI}X^{(n)}:=\coprod_{\mu} X^{\mu}. \]

Recall that as a graded vector space, the orbifold cohomology is defined as follows,
\[H^*_{\mathrm{orb}}(X^{(n)},\BQ):=H^{*-2\mathrm{age(\mu)}}(\CI X^{(n)},\BQ).\]
More explicitly, by (\ref{Inertia}),  we  get that 
\begin{equation} \label{orbifoldcoh}
	H^*_{\mathrm{orb}}(X^{(n)},\BQ)= \bigoplus_\mu H^{*-2\mathrm{age}(\mu)}(X^\mu,\BQ)^{C(\mu)},
\end{equation}
where $H^{*-2\mathrm{age}(\mu)}(X^\mu,\BQ)^{C(\mu)}$ denotes the $C(\mu)$-invariant part of the cohomology. 
\subsection{Invariants}	\label{secinv}
Let 
$ \vec{A}dm_{g,N}^{\epsilon}(X^{(n)},\beta)$ be the moduli space obtained from $Adm_{g,N}^{\epsilon}(X^{(n)},\beta)$ by putting a \textit{standard order}\footnote{We order the points in a fiber in accordance with the standard order of a ramification profile $\mu$.} on the fibers over marked points of the source curve.  The two moduli spaces are related as follows
\begin{equation} \label{decomposition1}
Adm_{g,N}^{\epsilon}(X^{(n)},\beta=	\coprod_{\underline \mu} [ \vec{A}dm^{\epsilon}_{g,N}(X^{(n)}, \beta,\underline \mu)/\prod_i \Aut(\mu^i)],
\end{equation}
where the union is taken over all $N$-tuples of partitions $(\mu^1, \hdots, \mu^N)$ of $n$. 
There exist 
evaluation morphisms at marked points

\begin{align*}\ev_{i}\colon   \vec{A}dm_{g,N}^{\epsilon}(X^{(n)},\beta) &\rightarrow \vec{\CI}X^{(n)}=\coprod_{\mu} X^{\mu}, \quad i=1, \dots, N, \\
f &\mapsto f_X(f_C^{-1}({x_i})).
\end{align*}
Consider the universal markings,
\[s_{i}\colon   \vec{A}dm_{g,N}^{\epsilon}(X^{(n)},\beta)\rightarrow \CC_{g,N,}\]
to the universal \textit{target} curve 
\[\CC_{g,N} \rightarrow   \vec{A}dm_{g,N}^{\epsilon}(X^{(n)},\beta).\]
We define the associated cotangent line bundles, 
\[\CL_{i}: =s^{*}_{i}(\omega_{\CC_{g,N} / Adm_{g,N}^{\epsilon}(X^{(n)},\beta)}), \quad i=1, \dots, N,\]
where $\omega_{\CC_{g,N} / Adm_{g,N}^{\epsilon}(X^{(n)},\beta)}$ is the universal relative dualising sheaf. We denote 
\[\psi_{i}:=\mathrm{c}_{1}(\CL_{i}).\]
With above structures at hand we can define $\epsilon$-admissible invariants. 
\begin{defn} The \textit{descendent} $\epsilon$-\textit{admissible invariants} are 
	\begin{multline*}\langle \psi_1^{m_{1}}\gamma_{1}, \dots, \psi_N^{m_{N}}\gamma_{N} \rangle^{\epsilon}_{g,\beta}:= \\
	 \frac{1}{\prod_i |\Aut(\mu^i)|}\int_{[ \vec{A}dm_{g,N}^{\epsilon}(X^{(n)},\beta)]^{\mathrm{vir}}}\prod^{i=N}_{i=1}\psi_{i}^{m_{i}} \ev^{*}_{i}(\gamma_{i}),
	 \end{multline*}
	where $\gamma_{1}, \dots, \gamma_{N} \in H^{*}(\vec{\CI} X^{(n)},\BQ)$ and $m_{1}, \dots m_{N}$ are non-negative integers. 
\end{defn}
\subsection{Relation to other invariants} \label{Relation2}
We will now explore how $\epsilon$-admissible invariants are related to the invariants associated to the spaces discussed in Section \ref{Relation1}.
\subsubsection{Classes} \label{Classes}
Let $\{\delta_1, \dots \delta_{m_X}\}$ be an ordered basis of $H^*(X,\BQ)$. Let  
\[\vec{\mu}=((\mu_1,\delta_{l_1}), \dots, (\mu_k,\delta_{l_k}))\]
be a cohomology-weighted partition of $n$. We put a standard order on $\vec{\mu}$, 
\[(\mu_{i}, \delta_{l_i}) > (\mu_{i'},\delta_{l_{i'}}),\]
if $\mu_{i}>\mu_{i'}$, or if $\mu_{i}=\mu_{i'}$ and $l_{i}>l_{i'}$. The underlying partition will be denoted by $\mu$. For each $\vec\mu$, we consider  a class 
\[\vec{\delta}=\delta_{l_1}\otimes\hdots \otimes \delta_{l_k} \in H^*(X^\mu,\BQ),\]
we then define 
\[\lambda(\vec{\mu}):= \frac{1}{\mathfrak{z}(\mu)} \pi_*(\delta_{l_1}\otimes\hdots \otimes \delta_{l_k})\in H_{\mathrm{orb}}^{*}(X^{(n)}, \BQ),\] 
where
\[ \pi\colon \vec{\CI} X^{(n)} \rightarrow \CI X^{(n)}\] 
is the natural projection. 
More explicitly, by using the identification 
\begin{equation} \label{identcoh}
H^*(X^{\mu},\BQ)^{C(\mu)} \cong H^*([X^{\mu}/C(\mu)],\BQ)\cong H^*(X^{\mu}/C(\mu),\BQ),
\end{equation}
the class $\lambda(\vec{\mu})$ is given by the following formula
\[ \frac{1}{\mathfrak{z}(\mu)} \sum_{h\in C(\mu)}h^*(\delta_{l_1}\otimes\hdots \otimes \delta_{l_k}) \in H^*(X^{\mu},\BQ)^{C(\mu)}.\] 
The importance of these classes is due to the fact they form a basis of $H_{\mathrm{orb}}^{*}(X^{(n)}, \BQ)$, see Proposition \ref{HilbSym}.
\subsubsection{Comparison} \label{comp} Given a collection of weighted partitions \[\vec \mu^i=((\mu^i_{1},\delta_{l^i_{1}}), \dots, (\mu^i_{k_{i}},\delta_{l^i_{k_i}})), \quad i=1,\dots, N,\]
where $\delta_{l^i_{j}} \in \{\delta_1, \dots \delta_{m_X}\}$.  The relative GW descendent invariants associated to the moduli space $\Mbar_{\mathsf{h}}(X\times C_{g,N}, (\gamma,n))$ agree with our invariants for $\epsilon=0^+$, 
\begin{multline*}
\langle \psi_1^{m_{1}}\vec \delta_1, \dots, \psi_N^{m_{N}}\vec \delta_N\rangle^{0^+}_{g,\beta} \\
=\frac{1}{\prod_i|\Aut(\mu^i)|}\ \int_{[\vec{M}_{\mathsf{m}}(X\times C_{g,N},(\gamma,n))]^{\mathrm{vir}}} \prod_{i=1}^N\psi_i^{m_i} \prod^{k_i}_{j=1} \ev_{i,j}^*(\delta_{l^i_j}),
\end{multline*}
such that we integrate over the space with a standard order on the ramification points, and 
\[\ev_{i,j}\colon  \vec{M}_{\mathsf{h}}(X\times C_{g,N},(\gamma,n)) \rightarrow X, \quad i=1, \dots, N, j=1, \dots, k_i,\] 
are evaluation morphisms defined by sending the corresponding point in a fiber over a marked point.

In the case of $\CK_{g,N}([X^{(n)}], \beta)$, we have evaluation morphisms to the inertia stack
\[\ev_i\colon \CK_{g,N}([X^{(n)}], \beta)\rightarrow \CI X^{(n)}, \quad i=1,\dots N,\]
for the definition of which we refer to \cite{AGV}.
We consider coarse $\psi$-classes on $\CK_{g,N}([X^{(n)}], \beta)$. Orbifold $\psi$-classes are rational multiples of coarse ones. Using the classes defined in the previous sections, the relevant invariants in this case are defined as follows, 

\[ \int_{[\CK_{g,N}([X^{(n)}], \beta)]^{\mathrm{vir}}} \prod_{i=1}^N \psi_i^{m_i} \ev_i^*(\lambda(\vec{\mu}^i)).\]
The next lemma concludes the comparison initiated in Section \ref{Relation1}.  
\begin{lemma} \label{invariantscomp} We have
		\[\langle \psi_1^{m_{1}}\vec \delta_1, \dots, \psi_N^{m_{N}}\vec \delta_N \rangle^{1}_{g,\beta}= \int_{[\CK_{g,N}([X^{(n)}], \beta)]^{\mathrm{vir}}} \prod_{i=1}^N \psi_i^{m_i} \ev_i^*(\lambda(\vec{\mu}^i)).\]

\end{lemma}
\textit{Proof.}
Consider the evaluation morphism to the coarse inertia stack,
\[ \ev_i \colon \CK_{g,N}([X^{(n)}],\beta) \rightarrow \CI X^{(n)} \rightarrow \coprod_\mu X^\mu/C(\mu).\]
By the identification (\ref{identcoh}), pulling back a class $\lambda(\vec{\mu}^i)$ from the coarse inertia stack is the same as pulling it back from the standard inertia stack. Moreover, the inertia stack parametrises representable maps from $B\BZ_r$ for different integers $r$, 
\[ B\BZ_r \rightarrow [X^{(n)}], \]
 see \cite[Section 3]{AGV}. By the universal property of $[X^{(n)}]$, a map $B\BZ_r \rightarrow [X^{(n)}]$ is given by a $S_n$-torsor on $B\BZ_r$ with a $S_n$-equivariant map to $X^n$. The evaluation morphism
\[
 \CK_{g,N}([X^{(n)}],\beta) \rightarrow \CI X^{(n)} 
 \]
 associates to the fiber of a $S_n$-torsor over a marked twisted point on a twisted curve, $B\BZ_r \in \CC$,  its image with respect to a map from that curve. Its image in the coarse quotient $X^\mu/C(\mu)$ after composing with 
 \[
 \CI X^{(n)} \rightarrow \coprod_\mu X^\mu/C(\mu)
 \]
is the fiber of the degree $n$ cover associated to the $S_n$-torsor via the correspondence discussed in Section \ref{Relation1}.
We therefore obtain that the evaluation morphisms commute, 
\[ 
\begin{tikzcd}[row sep=normal, column sep = normal]
	 \CK_{g,N}([X^{(n)}], \beta) \arrow{d}{\rho} \arrow{dr}{\mathrm{ev}_i} & \\
Adm_{g,N}^{1}(X^{(n)},\beta) \arrow{r}{\mathrm{ev}_i} & X^\mu/C(\mu)
\end{tikzcd}
\]
Finally, 
\[\pi^*(\lambda(\vec{\mu}))= \frac{1}{\mathfrak{z}(\mu)} \sum_{h\in C(\mu)}h^*(\delta_{l_1}\otimes\hdots \otimes \delta_{l_k}),\] 
where $\pi\colon \vec{\CI} X^{(n)} \rightarrow \CI X^{(n)}$; and the cardinality of $C(\mu)$ is exactly $\mathfrak z(\mu)$. Hence by pulling back classes from $Adm_{g,N}^{1}(X^{(n)},\beta)$  further to $\vec{A}dm_{g,N}^{1}(X^{(n)},\beta)$ and using that the integrals on the ordered space are invariant with respect to the permutation action of $C(\mu)$, we obtain the claim by Lemma \ref{fc}. 
   \qed

\section{Master space} \label{master}
\subsection{Definition of the master space}
The space $(0,1]$ of $\epsilon$-stabilities is divided into chambers, inside of which the moduli space $ Adm_{g,N}^{\epsilon}(X^{(n)},\beta)$ stays the same, and as $\epsilon$ crosses a wall between chambers, the moduli space changes discontinuously. By the definition of $\epsilon$-admissibility,  a wall is a number $1/d$ for some integer $d$. 

Let $\epsilon_0 \in (0,1]$ be a wall, and $\epsilon_+$,
$\epsilon_-$ be some values that are close to $\epsilon_0$ from the right and the left,  respectively. We set 
\[d_0=1/\epsilon_0\ \text{ and } \ \deg(\beta):=\mathsf {m}+\deg(\gamma)=d.\]

\begin{defn}
	A pre-admissible map $(P,C,\mathbf{x}, f)$ is called $\epsilon_0$\textit{-pre-admissible}, if 
	\begin{itemize}
		\item[(i)] for all points $p \in C$, 
		\[ \mathrm{mult}_p(\br(f))+\deg(f^*L_p)\leq 1/\epsilon_0,\]
		\item[(ii)] for all rational tails $T \subseteq C$, 
		\[\deg(\br(f)_{|T})+\deg(f^*L_{|T})\geq 1/\epsilon_0,\]
		\item[(iii)]  for all rational bridges $B \subseteq C$, the group of automorphisms of $f$ over $B$ are finite,\footnote{We consider automorphisms of $B$ together with its two special points; the same applies to components of the source curve over $B$. } 
		\[ 
		|\Aut(f_{|B})|<\infty;
		\]
	 and  the group of automorphisms on the source fixing $f$ is finite,  \[|\{(\phi \in \Aut(P)\mid f \circ \phi = f \}|<\infty.\]
		
	\end{itemize}
\end{defn}

Note that the last condition is saying that there are only finitely many automorphisms of the map with respect to the source and the target over rational bridges, and there are no unstable components contracted by the map $f$. We denote by $\mathfrak{Adm}^{\epsilon_0}_{g,N}(X^{(n)},\beta)$ the moduli space of $\epsilon_0$-pre-admissible maps with the specified discrete data. 

\begin{defn}
	Given a pre-admissible map $(P,C,f, \mathbf{x})$. We say that a rational tail $T \subseteq (C, \mathbf x)$ is of degree $d_0$, if 
	\[\deg(\br (f)_{|T})+\deg(f^*L_{|T})=d_0.\]
	We say a branching point $p \in C$ is of degree $d_0$, if 
	\[ \mathrm{mult}_p(\br(f))+\deg(f^*L_p) =d_0.\] 
\end{defn} 

Let $\FM^{ss}_{g,N,d}$ be the moduli space of weighted semistable curves of total degree $d$, i.e., curves with positive integers attached to each irreducible components whose sum is equal to $d$,  see  \cite[Definition 2.1.2]{YZ} for the precise definition.   There exists a map 
\[ \mathfrak{Adm}^{\epsilon_0}_{g,N}(X^{(n)},\beta) \rightarrow \FM^{ss}_{g,N,d}\] \[(P,C,f, \mathbf{x}) \mapsto (C,\mathbf{x},\underline{d}),\]
where the value of $\underline{d}$ on a subcurve $C' \subseteq C$ is defined as follows 
\[\underline{d}(C')=\deg(\br (f_{|C'}))+\deg(f^*L_{|C'}).\]

By $M\mathfrak{Adm}^{\epsilon_0}_{g,N}(X^{(n)},\beta)$, we denote the moduli space of $\epsilon_0$-pre-admissible maps with calibrated tails, defined as the fiber product
\[M\mathfrak{Adm}^{\epsilon_0}_{g,N}(X^{(n)},\beta)= \mathfrak{Adm}^{\epsilon_0}_{g,N}(X^{(n)},\beta)\times_{\FM^{ss}_{g,N,d}}M \widetilde{\FM}_{g,N,d},\] 
where $M \widetilde{\FM}_{g,N,d}$ is the moduli space of curves with calibrated tails introduced in \cite[Definition 2.8.2]{YZ}, which is a projective bundle  over a moduli space of curves with entangled tails, $\widetilde{\FM}_{g,N,d}$, see \cite[Section 2.2]{YZ}. The latter is constructed by induction on the integer $k$ as a sequence of blow-ups at the loci of curves with at least $k$ rational tails of degree $d_0$. 

A $B$-point in $M\mathfrak{Adm}^{\epsilon_0}_{g,N}(X^{(n)},\beta)$ is a tuple, 
\[(P,C, \mathbf{x},f, e, \CL, v_{1}, v_{2}),\]
where $(P,C, \mathbf{x},f,)$ is an $\epsilon_0$-pre-admissible map over $B$, and $(e, \CL, v_{1}, v_{2})$ is the calibration data. More explicitly, $e$ specifies a fiber of $\widetilde{\FM}_{g,N,d} \rightarrow \FM_{g,N,d}$ above the curve $[C] \in \FM_{g,N,d}(B)$, while $(\CL, v_{1}, v_{2})$ is a line bundle with sections $v_1 \in H^0(B, \BM_C \otimes \CL)$ and $v_2 \in H^0(B, \CL)$, where $\BM_C$ is the conormal bundle of the divisor of curves with at least one rational tail of degree $d_0$ (exceptionally in the case $(g,N)=(0,1)$, it is the cotangent line bundle associated to the marking). 

\begin{defn}
	We say that a rational tail $T \subseteq (C,\mathbf x)$ is \textit{constant}, if
	\[|\Aut(f_{|T})|=\Aut((P,C,f, \mathbf{x})_{|T})|=\infty.\]
\end{defn}
In other words, a rational tail $T \subseteq (C,\mathbf x)$ is constant, if at each connected component of $P_{|T}$, the map $f_{C|T}$ is equal to 
\[z^{\underline n}\colon (\sqcup^k \p^1)\cup P' \rightarrow \p^1\]
up to a change of coordinates, while $f_{X|T}$ is trivial. The notation is the same as in (\ref{constantmaps}). 
\begin{defn} A $B$-family of $\epsilon_0$-pre-admissible maps with calibrated tails
	\[(P,C, \mathbf{x},f, e, \CL, v_{1}, v_{2})\]
	is $\epsilon_0$-admissible, if 
	\begin{itemize}
		\item[(i)] every constant tail is an entangled tail,
		\item[(ii)] if a geometric fiber $C_b$ of $C$ over a point $b \in B$ has rational tails of degree $d_0$, then those rational tails contain all the degree-$d_0$ branching points,
		\item[(iii)] if $v_1(b)=0$, then $(P,C, \mathbf{x},f)_b$ is $\epsilon_+$-admissible,
		\item[(iv)] if $v_2(b)=0$, then $(P,C, \mathbf{x},f)_b$ is $\epsilon_-$-admissible. 
	\end{itemize} 
\end{defn}
Let 
\[M Adm^{\epsilon_0}_{g,N}(X^{(n)},\beta) \subset M\mathfrak{Adm}^{\epsilon_0}_{g,N}(X^{(n)},\beta)\]
 denote the moduli space of genus-$g$, $N$-marked, $\epsilon_0$-admissible maps with calibrated tails. 
\subsection{Obstruction theory}
The obstruction theory of $M Adm^{\epsilon_0}_{g,N}(X^{(n)},\beta)$  is constructed in the same way as the one of $ Adm^{\epsilon}_{g,N}(X^{(n)},\beta)$ in Section \ref{SectionObs}. More precisely,  the complex $R\pi_*(f_X^*T_X)$ defines a perfect obstruction theory relative to the stack $\FM_{\mathsf{h}}(\FC_{g,N}/ M\widetilde{\FM}_{g,N,d},n)$,
\[ \phi: R\pi_*(f_X^*T_X)^\vee \rightarrow \BL_{	MAdm_{g,N}^{\epsilon}(X^{(n)}, \beta))/\FM_{\mathsf{h}}(\FC_{g,N}/ M\widetilde{\FM}_{g,N,d},n)}.\]
Alternatively, the obstruction theory can also be defined via the moduli spaces of maps between twisted curve or maps between logarithmic curves.

\subsection{Properness} 

\begin{thm} A master space $M Adm^{\epsilon_0}_{g,N}(X^{(n)},\beta)$ is a quasi-separated Deligne--Mumford stack of finite type. 
\end{thm}
\textit{Proof.}  A master space $M Adm^{\epsilon_0}_{g,N}(X^{(n)},\beta)$ is simply a substack of the relative moduli space of maps to $\FC_{g,N}/ M\widetilde{\FM}_{g,N,d}$. The definition of $\epsilon_0$-admissibility ensures that it is Deligne--Mumford and of finite type, as in the proof of Proposition \ref{DM1}. See  \cite[Proposition 4.1.11]{YZ} for more details.
\qed
\\

We now deal with the properness of $M Adm^{\epsilon_0}_{g,N}(X^{(n)},\beta)$. We will follow the strategy of \cite[Section 5]{YZ}. Namely,  given a discrete valuation ring $R$ with the fraction field $K$. Let 
\[\xi^*=(P^*, C^*,\mathbf{x}^*,f^*, e^*, \CL^*, v^*_1,v^*_2) \in M Adm^{\epsilon_0}_{g,N}(X^{(n)},\beta)(K) \] be a family of $\epsilon_0$-admissible maps with calibrated tails over $\Spec K$. We will classify all the possible $\epsilon_0$-pre-admissible extensions of $\xi^*$ to $R$ up to a finite base change. There will be a unique one which is $\epsilon_0$-admissible. 

\subsubsection{The case of $(g,N,d)\neq(0,1,d_0)$} Assume $(g,N,d)\neq(0,1,d_0)$ and $\eta^*$ does not have rational tails of degree $d_0$.
Let 
\[\eta^*=(P^*, C^*, \mathbf{x}^*,f^*) \ \text{and} \ \lambda^*=(e^*, \CL^*, v^*_1,v^*_2)\]
be the underlying pre-admissible map and the calibration data of $\eta^*$, respectively. We define an extension of $\xi^*$ over $\Spec R$ (possibly after a finite base change), 
\[\xi_+=(\eta_+,\lambda_+) \in  M\FM^{\epsilon_0}_{g,N}(X^{(n)},\beta)(R),\]
as follows. Firstly,  the $\epsilon_0$-pre-admissible map
\[\eta_+= (P_+,C_+,\mathbf{x}_+,f_+),\] 
 is constructed in the same way as in \textit{Step 1,2} of the proof of Theorem \ref{properness} (possibly after a finite base change). More explicitly, we firstly apply the same construction as in \textit{Step 1} by separating the map $\eta^*$ into the contracted and non-contracted parts and taking the limit in $ \CK_{g,N}([X^{(n)}], \beta)$.  We then  apply modifications of \textit{Step 2} with respect to $\epsilon_+$-admissibility, leaving the degree-$d_0$ branching points which are limits of degree-$d_0$ branching points of the generic fiber untouched.  This means that we contract all rational tails which do not satisfy $\epsilon_+$-admissibility, which is possible by the assumption that the generic fiber is $\epsilon_0$-pre-admissible (i.e., a rational tail violating $\epsilon_+$ admissibility cannot be a limit of a rational tail at the generic fiber), and blow up whenever there is an excess contribution from contracted components and nodes not satisfying $\epsilon_+$-admissibility except if this is a degree-$d_0$ branching point coming from the generic fiber. 

The family $\eta_+$ is the one closest to being $\epsilon_+$-admissible limit of $\eta^*$ (i.e., the only obstacle to being $\epsilon_+$-admissible is the degree-$d_0$ branching points in the generic fiber). The calibration $\lambda_+$ is given by the unique
extension of $\lambda^*$ to the curve $C_+$, which exists by \cite[Lemma 5.1.1 (1)]{YZ}. 

Let 
\[\{p_i \mid  i= 1, \dots, \ell \}\]
be an ordered set, consisting of nodes of degree-$d_0$ rational tails and degree-$d_0$ branching points of the central fiber 
\[p_i\in C_{+|\Spec \BC}\subset C_{+}.\] 
We now define 
\[b_i \in \BR_{>0} \cup \{\infty\}, \ i=1, \dots, \ell  \] 
as follows.
Set $b_i$
to be $\infty$, if $p_i$ is a degree-$d_0$ branching point. If $p_i$ is a node of a rational tail, then we define $b_i$ via the singularity type of $C_+$ at $p_i$. Namely, if the total space of the family $C_+$ has an $A_{b-1}$-type singularity at $p_i$, we set $b_i=b$.  \\

We now classify all $\epsilon_0$-pre-admissible modifications of $\xi_+$ in the sense of Definition \ref{modifiation}. By \cite[Lemma 5.1.1 (1)]{YZ}, it is enough to classify the modifications of $\eta_+$, because the calibration data extends uniquely.  

All modifications of $\eta_+$ are given by blow-ups and contractions around the points $p_i$ after taking base-changes with respect to finite extensions of $R$. The result of these modifications will be a change of singularity type of $\eta_+$ around $p_i$. Hence the classification will depend on an array of rational numbers 
\[\underline{\alpha}=(\alpha_1,\dots, \alpha_\ell)\in \BQ_{\geq 0}^\ell,\] the nominator of which keeps track of the singularity type around $p_i$, while the denominator is responsible for the degree of an extension of $R$. The precise statement is the following lemma. 
\begin{lemma} \label{neq01d}
	For each $\underline{\alpha}=(\alpha_1,\dots, \alpha_\ell)\in \BQ_{\geq 0}^\ell$, such that $\underline{\alpha}\leq \underline{b}$, there exists an $\epsilon_0$-pre-admissible modification $\eta_{\underline{\alpha}}$ of $\eta_+$ 
	with the following properties.
	\begin{itemize}
		\item Up to a finite base change, 
		\[\eta_{\underline{\alpha}}\cong \eta_{\underline{\alpha}'} \iff \underline{\alpha}=\underline{\alpha}' .\]
		\item  Given an $\epsilon_0$-pre-admissible modification $\tilde{\eta}$ of $\eta_+$, then, up to a finite base change,  there exists $\underline{\alpha}$ such that  \[\tilde{\eta} \cong \eta_{\underline{\alpha}}.\] 
		\item  the central fiber of $\eta_{\underline{\alpha}}$ is $\epsilon_+$-admissible, if and only if $\underline{a}=\underline{b}$. 
	\end{itemize}
\end{lemma}

\textit{Proof.}
Let us choose a fractional presentation of $(a_1,\dots, a_\ell)$ with a common denominator
\[(a_1,\dots, a_\ell)=(\frac{a'_1}{r},\dots, \frac{a'_\ell}{r}).\]
Take a base change of $\eta_+$ with respect to a degree-$r$ extension of $R$. We then construct $\eta_{\underline{\alpha}}$ by applying modifications to $\eta_+$ around each point $p_i$, the result of which is a family
\[\eta_{\alpha_i}=(P_{\alpha_i}, C_{\alpha_i}, \mathbf{x}_{\alpha_i},f_{\alpha_i}),\]
constructed as follows. 
\\

\textit{Case 1.} If $p_i$ is a node of a  degree-$d_0$ rational tail and $a_i\neq 0$, we blow up $C_+$ at $p_i$,
\[\mathrm{Bl}_{p_i}(C_+) \rightarrow C_+.\]
The map $f_{C_+}$ then defines a rational map 
\[f_{C_+}\colon P_+ \dashrightarrow \mathrm{Bl}_{p_i}(C_+).\]
We can eliminate the indeterminacies of the map above by blowing up $P_+$ to obtain an everywhere-defined map 
\[f_{\mathrm{Bl}_{p_i}(C_+)} \colon \widetilde{P}_+ \rightarrow \mathrm{Bl}_{p_i}(C_+),\] 
we take a minimal blow-up with such property.
The exceptional curve $E$ of $\mathrm{Bl}_{p_i}(C_+)$ is a chain of $rb_i$ rational curves. The exceptional curve of $\widetilde{P}_+$ is a disjoint union $\sqcup E_j$, where each $E_j$ is a chain of $r b_i$ rational curves mapping to $E$ without contracted components; on each irreducible component of $E_j$ the map is given by $z^k \colon \p^1 \rightarrow \p^1$, where $k$ is determined by the ramification over the node before the blow-up.  We contract all the rational curves but the $a'_i$-th ones in both $E$ and $E_j$ for all $j$. The resulting families are $C_{\alpha_i}$ and $P_{\alpha_i}$, respectively. The family $C_{\alpha_i}$ has an $A_{\alpha_i'-1}$-type singularity at $p_i$. The marking $\mathbf{x}_+$ clearly extends to a marking $\mathbf{x}_{\alpha_i}$ of $C_{\alpha_i}$. The map $f_{\mathrm{Bl}_{p_i}(C_+)}$ descends to a map
\[f_{C_{\alpha_i}}\colon P_{\alpha_i} \rightarrow C_{\alpha_i}.\] 
The map $f_{+,X}$ is carried along with all those modifications to a map 
\[f_{\alpha_i,X} \colon P_{\alpha_i} \rightarrow X,\] because exceptional divisors are of degree $0$ with respect to $f_{+,X}$, hence the contraction of curves in the exceptional divisors does not introduce any indeterminacies. We thereby constructed the family $\eta_{\alpha_i}$. 
\\

\textit{Case 2.} Assume now that $p_i$ is a node of a degree-$d_0$ rational tail, but $a_i=0$. The family $C_{\alpha_i}$ is then given by the contraction of that degree-$d_0$ rational tail, which is possible because it is not a limit of a degree-$d_0$ rational tail at the generic fiber. The resulting family $C_{\alpha_i}$ is smooth at $p_i$. The marking $\mathbf{x}_+$ extends to a marking $\mathbf{x}_{\alpha_i}$ of $C_{\alpha_i}$. The family $P_{\alpha_i}$ is set to the stabilisation of  $P_+$ together with the map (i.e., we contract rational components which become unstable after the contraction of the rational tail in the target). The map $f_{\alpha_i}$ is the composition of the contraction and $f_+$. 
\\

\textit{Case 3.} If $p_i$ is a branching point, we blow up $C_+$ inductively $a'_i$ times, starting with a blow-up at $p_i$ and then continuing with a blow-up at a point of the exceptional curve of the previous blow-up. We then contract all rational curves in the exceptional divisor but the last one. The resulting family is $C_{\alpha_i}$, it has an $A_{a'_i-1}$-type singularity at $p_i$. The marking $\mathbf{x}_+$ extends to the marking $\mathbf{x}_{\alpha_i}$ of $C_{\alpha_i}$. The map $f_{C_+}$ then defines a rational map
\[f_{C_+} \colon P_+ \dashrightarrow C_{\alpha_i}.\]
We set 
\[f_{C_{\alpha_i}} \colon P_{\alpha_i} \rightarrow C_{\alpha_i}\]
to be the minimal resolution of indeterminacies of the rational map above. More specifically, $P_{\alpha_i}$ is obtained by consequently blowing up $P_+$ and contracting all the rational curves in the exceptional divisors but the last ones. The map $f_{+,X}$ is carried along, as in \textit{Case 1}. 
\\

After having applied all modifications associated to $\underline{\alpha}\in \BQ_{\geq 0}^\ell$ to points $p_i$, we obtain the family $\eta_{\underline{\alpha}}$. It is indeed $\epsilon_0$-pre-admissible, since all modifications take place at degree-$d_0$ rational tails or degree-$d_0$ branching points by either changing the singularity type of their attachment or moving a degree-$d_0$ point to a degree-$d_0$ rational tail. Up to a finite base change, the resulting family is uniquely determined by $\underline{\alpha}=(\alpha_1,\dots, \alpha_\ell)$ and independent of its fractional presentation, because of the singularity types at points $p_i$ and the degree of an extension of $R$.  

Given now an arbitrary $\epsilon_0$-pre-admissible modification 
\[\eta=(P,C,\mathbf{x},f)\] of $\eta_+$. Possibly after a finite base change, there exists a modification 
\[\tilde \eta=(\widetilde{P}, \widetilde{C}, \tilde{\mathbf{x}},\tilde{f})\] that dominates both $\eta$ and $\eta_+$ in the sense of $(\ref{dominatin})$, and which is constructed in the same way. We take a minimal modification with such property.   By the assumption of minimality and $\epsilon_0$-pre-admissibility of $\eta$ and $\eta_+$, the families $(\widetilde{C},\widetilde{P})$ and $(C_+,P_+)$ must be related by blow-ups over $p_i$,
\[ \widetilde C \rightarrow  C_+ \ \text{ and } \ \widetilde P \rightarrow  P_+, \]
 since all other blow-ups would introduce either unstable rational bridges or rational tails of wrong degree contradicting the $\epsilon_0$-pre-admissibility of $\eta$. Moreover, again by the $\epsilon_0$-pre-admissibility of $\eta$, the projections 
\[\widetilde C \rightarrow  C \ \text{ and } \ \widetilde P \rightarrow  P\]
are given by contraction of degree-$d_0$ rational tails or rational components which do not satisfy $\epsilon_0$-pre-admissibility. These are exactly the operations described in \textit{Steps 1,2,3} of the proof, implying that $(C,P)$ must be one of the families constructed previously.  Uniqueness of maps follows from the seperatedness of the moduli space of maps to a fixed target. Hence we obtain that 
\[\eta \cong \eta_{\underline{\alpha}}\]
for some $\underline{\alpha}=(\alpha_1,\dots, \alpha_\ell)\in \BQ_{\geq 0}^\ell$, where $\underline{\alpha}$ is determined by the singularity types of $\eta$ at points $p_i$. 

The last claim follows from the construction of $\eta_+$ and the fact that we set $b=\infty$ for all degree-$d_0$ branching points of the central fiber. In particular, if $\alpha_i=b$ for all $i$, where $\alpha_i \in \BQ_{\geq 0}$, this means that we do not modify $\eta_+$ and $\eta_+$ does not have  degree-$d_0$ branching points, hence by construction it is $\epsilon_+$-admissible. 
\qed

\subsubsection{The case of $(g,N,d)=(0,1,d_0)$} We now assume that $(g,N,d)=(0,1,d_0)$. In this case, the calibration bundle is the relative cotangent bundle along the unique marking, see \cite[Definition 2.8.1]{YZ}. Moreover, there is no entanglement. Given a family of pre-admissible maps $(P,C,\mathbf{x}, f)$, we will denote the calibration bundle by $\BM_{C}$.   The calibration data $\lambda$ is given by a rational section $s$ of $\BM_{C}$. 

Let 
\[ \xi_+=(\eta_+,\lambda_+) \in  M\mathfrak{Adm}^{\epsilon_0}_{0,1}(X^{(n)},\beta)(R)\]
be the family over  $\Spec R$ (possibly after a finite base change), such that $\eta_+$ is an $\epsilon_+$-admissible extension of $\eta^*$, if there is no degree-$d_0$ branching point in $\eta^*$; it is constructed  as in the case  $(g,N,d)\neq(0,1,d_0)$. If there is a degree-$d_0$ branching point, then the map $\eta^*$ is a contraction of a curve to a point in $\p^1$.  In this case, we let $\eta_+$ be any pre-admissible extension (i.e., we just need to take the limit of that point to which the curve is contracted). The calibration data $\lambda_+$ is given by a rational section $s_+$ which is an extension of the section $s^*$ of $\BM_{C^*}$ to $\BM_{C_+}$. In general, given a modification $\widetilde{\eta}$ of $\eta_+$ over a degree-$r$ extension of $R$, the section $s^*$ extends to a rational section $\tilde{s}$ of $\BM_{\widetilde{C}}$. 
\begin{defn}
	The order of a modification $\widetilde{\eta}$ of $\eta_+$ is defined to be $\mathrm{ord}(\tilde{s})/r$ at the closed point of $\Spec R$.
\end{defn}
We set $b=\mathrm{ord}(s_+)$, if there is no degree-$d_0$ branching point in the generic fiber of $\eta^*$. Otherwise set $b=-\infty$.

\begin{lemma} \label{01d}
	For each $\alpha \in \BQ$, such that $\alpha \geq b$, there exists an $\epsilon_0$-pre-admissible modification $\eta_{\alpha}$ of $\eta_+$ of order $\alpha$ 
	with the following properties.
	\begin{itemize}
		\item Up to a finite base change, 
		\[\eta_{\alpha}\cong \eta_{\alpha'} \iff \alpha=\alpha'.\]
		\item  Given a $\epsilon_0$-pre-admissible modification $\tilde{\eta}$ of $\eta_+$, then, up to a finite base change,  there exists $\alpha$ such that  \[\tilde{\eta} \cong \eta_{\alpha}.\]
		
		\item  The central fiber of $\eta_{\alpha}$ is $\epsilon_+$-admissible, if and only if $\alpha=b$. 
	\end{itemize} 
\end{lemma}


\textit{Proof.} Assume  $\eta^*$ does not have a degree-$d_0$ branching point. We choose a fractional presentation 
\[
\alpha=\frac{\alpha'}{r}.
\]  
We take a base change of $\eta_+$ with respect to a degree-$r$ extension of $R$. We  blow up consequently $\alpha'$ times the central fiber of $C_+$ at the unique marking. We then contract all rational curves in the exceptional divisor but the last one. The resulting family with markings is $(C_{\alpha},\mathbf{x}_{\alpha})$. We do the same with the family $P_+$ at the points in the fiber over the marked point to a obtain the family $P_{\alpha}$ and the map 
\[f_{P_{\alpha}}\colon P_{\alpha} \rightarrow C_{\alpha},\]
the map $f_{+,X}$ is carried along with blow-ups and contractions. The resulting family of $\epsilon_0$-pre-admissible maps is of order $\alpha$. 

Assume now that the generic fiber has a degree-$d_0$ branching point. We take the base change of $\eta_+$ with respect to a degree-$r$ extension of $R$. After choosing some trivialisation of $\BM_{C^*}$, we have that \[s^*=\pi^{r \alpha_+} \in K,\]
where $\alpha_+$ is the order of vanishing of $s_+$ before the base-change and $\pi$ is some uniformiser of $R$. Because of automorphisms of $\p^1$ which fix a branching point and a marked point, we have an isomorphisms of $\epsilon_0$-pre-admissible maps with calibrated tails,
\[(\eta^*,s^*)\cong (\eta^*, \pi^{c} \cdot s^*) \]
for arbitrary $c\in \BZ$. Hence we can multiply the section $s_+$ with $\pi^{\alpha'-r\alpha_+}$ to obtain a modification of order $\alpha$. 

The facts that these modifications classify all possible modifications of $\eta_+$ and $\eta_\alpha$ is $\epsilon_+$-admissible for $\alpha=b$ follow from the same arguments as in the case $(g,N,d)\neq (0,1,d_0)$. 
\qed

\begin{thm} \label{Masterproper}  The moduli space $M Adm^{\epsilon_0}_{g,N}(X^{(n)},\beta)$ is proper.
\end{thm}
\textit{Proof.} First, using classifications of modifications of $\eta_+$ of Lemma \ref{neq01d}, \ref{01d}, extensions of $\eta^*$ with degree-$d_0$ rational tails for $(g,N,d)\neq(0,1,d_0)$ can also be classified using  arguments from \cite[Section 5.4]{YZ}. The classification takes exactly the same form with the difference that we have to separate $\eta^*$ into the part without degree-$d_0$ rational tails and the part consisting of  degree-$d_0$ rational tails. Then using both classifications discussed previously, extensions of such $\eta^*$ can be classified in the same way by means of a vector of rational numbers $\underline{\alpha}$. 

 The proof of properness of the master space then follows  \cite[Proposition 5.0.1]{YZ}, and uses the classifications of modifications of $\eta_+$ from Lemma \ref{neq01d}, \ref{01d}, as well as for the case when  $(g,N,d)\neq(0,1,d_0)$ and $\eta^*$  has a degree-$d_0$ rational tail. More explicitly, an $\epsilon_0$-pre-admissible modification of $\eta_+$ is $\epsilon_0$-admissible, if and only if it satisfies the inequalities from \cite[Equation (5.2), (5.3)]{YZ}. This shows both the existence and uniqueness parts of the valuative criteria of properness.     \qed

\section{Wall-crossing} \label{masterSym}
\subsection{Definition of $I$-functions} \label{graphspaceSym}
For a class $\beta = (\gamma,\mathsf m)\in H_2(X,\BZ)\oplus \BZ$, consider now 
\[\Mbar_{\mathsf{m}}(X\times \p^1/X_{\infty}, (\gamma,n)),\]
which is the space of relative stable maps with possibly disconnected sources of degree $(\gamma,n)$ to $X\times \p^1$ relative to 
\[X_{\infty}:=X\times \{\infty \} \subset X\times \p^1,\]
such that $\infty=[1,0] \in \p^1=\p(\BC^2)$. It admits a decomposition into subspaces of maps with the ramification profile $\mu$ over $\infty$, 

\[\Mbar_{\mathsf{m}}(X\times \p^1/X_{\infty}, (\gamma,n)= \coprod_{\mu}\Mbar_{\mathsf{m}}(X\times \p^1/X_{\infty}, (\gamma,n),\mu).  \]
Note that we fix the degree of the branching divisor $\mathsf m$ instead of the genus $\mathsf h$, the two determine each other by Lemma \ref{RHformula}. 

Consider the $\BC^*$-action on $\p^1$ given by 
\[t[x,y]=[tx,y], \ t\in \BC^*.\]
We denote the class of the weight 1 representation $\BC_{\mathrm{std}}$ of $\BC^*$ as follows,
\[ 
e_{\BC^*}(\BC_{\mathrm{std}})=z. 
\]
The $\BC^*$-action on $\p^1$ induces a $\BC^*$-action on $\Mbar_{\mathsf m}(X\times \p^1/X_{\infty}, (\gamma,n))$. Let 
\[F_{\beta} \subset \Mbar_{\mathsf{m}}(X\times \p^1/X_{\infty}, (\gamma,n))^{\BC^*}\] 
be the distinguished $\BC^*$-fixed component consisting of maps to $X\times \p^1$ (no expanded degenerations).  Said differently, $F_{\beta}$ is the moduli space of maps, which are admissible over $\infty \in \p^1$ and whose degree lies entirely over $0 \in \p^1$ in the form of a branching point. 

The virtual fundamental class of $F_{\beta}$,
\[[F_{\beta}]^{\mathrm{vir}} \in H_*(F_{\beta},\BQ),\] 
is defined via the fixed part of the perfect obstruction theory of 
\[\Mbar_{\mathsf{m}}(X\times \p^1/X_{\infty}, (\gamma,n)).\] The virtual normal bundle $N_{F_{\beta}}^{\mathrm{vir}}$ is defined via the moving part of the obstruction theory. By putting a standard order on the ramification points over $\infty$,  we obtain an ordered space $\vec{F}_{\beta}$, and the associated evaluation morphism
\[\mathsf{ev} \colon \vec{F}_{\beta} \rightarrow \vec{\CI} X^{(n)}\]
defined in the same way as in Section \ref{secinv}. 
\begin{defn}
We define a class in $H^{*}(\vec{\CI} X^{(n)},\BQ)[z^{\pm}]\otimes_{\BQ}\BQ[\![q^{\beta}]\!]$, called $I$-function,
\[I(z)=1+ \prod_j \mu_j\cdot  \sum_{\beta\neq0}q^{\beta}  \evsf_{*} \left(\frac{[\vec{F}_{\beta}]}{e_{\BC^*}(N_{\vec{F}_{\beta}}^{\mathrm{vir}})}\right),\]
where  $\prod_j\mu_j$ is a locally constant function on connected components of $\vec{\CI} X^{(n)}$ labelled by partitions $\mu$; we expand rational functions  in $z$ appearing in the expression above in the range $|z|>1$.
Let 
\[I_+(z) \in H^{*}(\vec{\CI} X^{(n)},\BQ)[z]\otimes_{\BQ}\BQ[\![q^{\beta}]\!]\] 
be the truncation $[zI(z)-z]_+$ by taking only non-negative powers of $z$. Let 
\[I_{+,\beta}(z) \in H^{*}(\vec{\CI} X^{(n)},\BQ)[z]\]
be the coefficient of $I_+(z)$ at $q^{\beta}$. For later, it is also convenient to define the class
\[\CI_{\beta}:=\frac{\prod_j{\mu_j} }{e_{\BC^*}(N^\mathrm{vir}_{\vec{F}_{\beta}})} \in H^*(\vec{F}_{\beta},\BQ)[z^{\pm}].\]
\end{defn}

\subsection{Wall-crossing formula} \label{secwall} From now on, we assume that 
\[2g-2+n+1/d_0\cdot \deg(\beta)>0,\]
for the case of $(g,N,d_0)=(0,1,d_0)$ we refer to \cite[Section 6.4]{YZ}. There exists a natural $\BC^*$-action on the master space $M Adm^{\epsilon_0}_{g,N}(X^{(n)},\beta)$ given by  
\[t \cdot (P, C,\mathbf{x},f, e, \CL, v_1,v_2)= (P, C,\mathbf{x},f, e, \CL, t \cdot v_1,v_2), \quad t \in \BC^*.\] 
By arguments presented in \cite[Section 6]{YZ}, there are three types of $\BC^*$-fixed $\epsilon_0$-admissible maps with calibrated tails: 
\begin{enumerate}
	\item $v_2=0$ and $(P, C,\mathbf{x},f)$ is  $\epsilon_-$-admissible,
	\item $v_1=0$ and $(P, C,\mathbf{x},f)$ is $\epsilon_+$-admissible,
	\item $v_1\neq 0$ and $v_2\neq0 $, and all degree-$d_0$ rational tails of $(P, C,\mathbf{x},f, e)$ are entangled constant rational tails. 
\end{enumerate}
The $\BC^*$-fixed locus therefore admits the following expression  \[M Adm^{\epsilon_0}_{g,N}(X^{(n)},\beta)^{\BC^*}=F_- \sqcup F_+ \sqcup \coprod_{\underline{\beta}} F_{\underline{\beta}},\] 
where each of the components correspond to one of the types of $\BC^*$-fixed maps described above. We will now explain the precise meaning of each term in the union above, giving a description of virtual fundamental classes and virtual normal bundles. 
\subsubsection{$F_-$} This is the simplest component, 
\[
F_-= Adm^{\epsilon_-}_{g,N}(X^{(n)},\beta), \quad
N_{F_-}^{\mathrm{vir}}=\BM^{\vee}_-,\]
where $\BM^{\vee}_-$ is the dual of the calibration bundle $\BM_-$ on $ Adm^{\epsilon_-}_{g,N}(X^{(n)},\beta)$, with a trivial $\BC^*$-action of weight $-1$, cf. \cite[Section 6.2]{YZ}. The obstruction theories also match, therefore 
\[[F_-]^{\mathrm{vir}}=[ Adm^{\epsilon_-}_{g,N}(X^{(n)},\beta)]^{\mathrm{vir}}\]
with respect to the identification above. 
\subsubsection{$F_+$} We define
\[\widetilde{A}dm^{\epsilon_+}_{g,N}(X^{(n)},\beta):=  Adm^{\epsilon_+}_{g,N}(X^{(n)},\beta) \times_{\FM_{g,N,d}} \widetilde{\FM}_{g,N,d},\] 
then
\[
F_+=\widetilde{A}dm^{\epsilon_+}_{g,N}(X^{(n)},\beta),\quad
N_{F_+}^{\mathrm{vir}}=\BM_+ ,
\]
where, as previously, $\BM_+$ is the calibration bundle on $\widetilde{A}dm^{\epsilon_+}_{g,N}(X^{(n)},\beta)$ with the trivial $\BC^*$-action of weight 1. The obstruction theories also match and 
\[p_*[\widetilde{A}dm^{\epsilon_+}_{g,N}(X^{(n)},\beta)]^{\mathrm{vir}}=[ Adm^{\epsilon_+}_{g,N}(X^{(n)},\beta)]^{\mathrm{vir}},\]
where \[p \colon \widetilde{A}dm^{\epsilon_+}_{g,N}(X^{(n)},\beta) \rightarrow  Adm^{\epsilon_+}_{g,N}(X^{(n)},\beta) \] 
is the natural projection. This follows from \cite[Theorem 5.0.1]{Co} and the fact that $\widetilde{\FM}_{g,N,d}$ is a blow-up of $\FM_{g,N,d}$ in a smooth locus, and the perfect obstruction theory of $\widetilde{A}dm^{\epsilon_+}_{g,N}(X^{(n)},\beta)$ is defined relatively to $\widetilde{\FM}_{g,N,d}$, 
\subsubsection{$F_{\underline{\beta}}$} These are the wall-crossing components, which will be responsible for wall-crossing formulas. Let 
\[\underline{\beta}=(\beta',\beta_1,\dots,\beta_k)\] be a $k+1$-tuple of classes in $H_2(X,\BZ) \oplus \BZ$, such that $\beta=\beta'+\beta_1+ \dots + \beta_k$ and $\deg(\beta_i)=d_0$, where the degree of a class $\beta$ is defined in the beginning of Section \ref{master}.   Then a component $F_{\underline{\beta}}$ is defined as follows
\begin{align*}
F_{\underline{\beta}}=\{ \xi \mid \ &\xi \text{ has exactly }k \text{ entangled tails,} \\
&\text{which are constant and of degree } \beta_1, \dots, \beta_k \}.
\end{align*}
Let 
\[
\begin{tikzcd}[row sep=small, column sep = small]
\CE_i \arrow[r] &  F_{\underline{\beta}} \arrow[l, bend left=20,"p_i"] &i=1,\dots, k,
\end{tikzcd}
\]
be the universal $i$-th entangled rational tail with the universal marking $p_i$  given by the node. We define $\psi(\CE_i)$ to be the $\psi$-class associated to the marking $p_i$.  Let 
\[\widetilde{\mathrm{gl}}_k \colon \widetilde{\FM}_{g,N+k, d-kd_0}\times (\FM_{0,1,d_0})^k \rightarrow \widetilde{\FM}_{g,N,d}\]
be the gluing morphism, cf. \cite[Section 2.4]{YZ}.
Let 
\[\mathfrak{D}_i \subset \widetilde{\FM}_{g,N,d}\]
be the divisor defined as the closure of the locus of curves with exactly $i+1$ entangled tails. 
Finally, let 
\[Y \rightarrow \widetilde{A}dm^{\epsilon_+}_{g,N}(X^{(n)},\beta')\] be the stack of $k$-roots of $\BM_+$. In the following proposition, we use spaces with ordered ramification points over $N$ marked points, which acquire an arrow superscript.

\begin{prop} \label{isomophism}
There exists a canonical isomorphism,
\[\widetilde{\mathrm{gl}}_k^*F_{\underline{\beta}} \cong \vec{Y} \times_{(\vec{\CI}X^{(n)})^k} \prod_{i=1}^{k} \vec{F}_{ \beta_i}/\prod^{k}_{i=1} \Aut(\mu^{N+i}).\]
With respect to the isomorphism above, we have 
\begin{align*}[\widetilde{\mathrm{gl}}_k^* F_{\underline{\beta}}]^{\mathrm{vir}}=&  [\vec{Y}]^{\mathrm{vir}} \times_{(\vec{\CI}X^{(n)})^k} \prod_{i=1}^{k} [\vec{F}_{ \beta_i}]^{\mathrm{vir}}/\prod^{k}_{i=1}|\Aut(\mu^{N+i})|, \\	\frac{1}{e_{\BC^*}(\widetilde{\mathrm{gl}}_k^* N_{F_{\underline{\beta}}}^{\mathrm{vir}})}=&\frac{\prod^k_{i=1}(z/k+\psi(\CE_i) )}{-z/k-\psi(\CE_1)-\psi_{n+1}-\sum^{\infty}_{i=k}\mathfrak{D}_i} \cdot \prod^k_{i=1}  \CI_{\beta_i}(z/k+\psi(\CE_i)).
\end{align*}
\end{prop}
\textit{Proof.} The first claim follows from the same arguments as in  \cite[Lemma 6.5.5]{YZ}, and the description of the $\BC^*$-fixed components $F_{\underline{\beta}}$. In our case, to glue maps, we also need to order the ramification points of the source curves, this gives rise to ordered spaces in the statement of the proposition.    We quotient by $\Aut(\mu^i)$ to remove the order on the ramification points after gluing is done.

To establish the second claim, we need to analyse the fixed and moving parts of the obstruction theory.  Following \cite[Lemma 6.5.5]{YZ}, we split the obstruction theory of the master space into the obstruction theory of maps to a fixed target curve and the obstruction theory of target curves, 
\begin{equation*}
	\BT^{\mathrm{vir}}_{M Adm^{\epsilon_0}_{g,N}(X^{(n)},\beta)/M\widetilde{\FM}_{g,N,d},}  \rightarrow \BT^{\mathrm{vir}}_{M Adm^{\epsilon_0}_{g,N}(X^{(n)},\beta)} \rightarrow \BT_{M\widetilde{\FM}_{g,N,d}} \rightarrow.
\end{equation*}
The contribution of  $\BT_{M\widetilde{\FM}_{g,N,d}}$ is independent of the setting (let it be quasimaps or our admissible maps), its moving part leads to the term \[
\frac{\prod^k_{i=1}(z/k+\psi(\CE_i) )}{-z/k-\psi(\CE_1)-\psi_{n+1}-\sum^{\infty}_{i=k}\mathfrak{D}_i}
\] in the expression from the statement of the proposition. The contribution  of
\[
\BT^{\mathrm{vir}}_{M Adm^{\epsilon_0}_{g,N}(X^{(n)},\beta)/M\widetilde{\FM}_{g,N,d},}
\]
 is responsible for the rest, and the only difference with the setting of \cite{YZ} is that the spaces of quasimaps from $\p^1$ are replaced by spaces of maps to $X \times \p^1$, which play the completely analogues role in our setting. In particular, the factors  \[
 \CI_{\beta_i}(z/k+\psi(\CE_i))
 \]
 arise as its moving part.

  The factor $ \prod_j{\mu_j^i}$ in the definition of $\CI_{\beta_i}(z/k+\psi(\CE_i))$ appears in the localisation formula due to the splitting of nodes of curves, e.g.,  see \cite[Theorem 3.6]{GV}, the same analysis holds in our setting. 
\qed
\begin{thm} \label{wallcrossingSym} Assuming $2g-2+N+1/d_0\cdot\deg(\beta)>0$, we have 
\begin{multline*}
	\langle \psi_1^{m_{1}}\gamma_{1}, \dots, \psi^{m_{N}}_N\gamma_{N} \rangle^{\epsilon_{-}}_{g,\beta}-\langle \psi_1^{m_{1}}\gamma_{1}, \dots, \psi^{m_{N}}_N\gamma_{N} \rangle^{\epsilon_{+}}_{g, \beta}\\
	=\sum_{k\geq1}\sum_{\underline{\beta}} \frac{1}{k!}\langle \psi_1^{m_{1}}\gamma_{1}, \dots, \psi^{m_{N}}_N\gamma_{N}, I_{+,\beta_1}(-\psi_{N+1}), \hdots, I_{+,\beta_k}(-\psi_{N+k}) \rangle^{\epsilon^+}_{g, \beta'},
\end{multline*}
where $\underline{\beta}$ runs through all the $(k+1)$-tuples of effective curve classes 
\[\underline{\beta}=(\beta', \beta_{1}, \dots, \beta_{k}),\]
such that $\beta=\beta'+ \beta_{1}+ \dots + \beta_{k}$ and $\deg(\beta_{i})=d_{0}$ for all $i=1, \dots, k$.
\end{thm}
\textit{Sketch of Proof.} 
We will explain the master-space technique, which is simply a localisation formula on the master space.  For all details we refer to \cite[Section 6]{YZ}, the shape of the localisation formula in our case is completely identical. 

First, for brevity, let 
\[\alpha=\prod^{i=N}_{i=1}\psi_{i}^{m_{i}}\ev^{*}_{i}(\gamma_{i})\in H^*(M \vec{A}dm^{\epsilon_0}_{g,N}(X^{(n)},\beta),\BQ)\] 
be the class corresponding to descendent insertions. By the virtual localisation formula, we obtain 
\begin{multline*} [M \vec{A}dm^{\epsilon_0}_{g,N}(X^{(n)},\beta)]^{\mathrm{vir}} \cap \alpha \\
=  \sum \iota_{\vec{F}_\star *}\left( \frac{[\vec{F}_\star]^{\mathrm{vir}}\cap \alpha}{e_{\BC^*}(N_{\vec{F}_\star}^{\mathrm{vir}})} \right) \in H^{\BC^*}_*(M \vec{A}dm^{\epsilon_0}_{g,N}(X^{(n)},\beta)[z^{-1}],
\end{multline*}
where $\vec{F}_{\star}$ are the components of the $\BC^*$-fixed locus of $M \vec{A}dm^{\epsilon_0}_{g,N}(X^{(n)},\beta)$ with ordered ramification points above $N$ markings on target curves.  Since the equivariant virtual fundamental class $[M \vec{A}dm^{\epsilon_0}_{g,N}(X^{(n)},\beta)]^{\mathrm{vir}}$ is contained in $H^{\BC^*}_*(M \vec{A}dm^{\epsilon_0}_{g,N}(X^{(n)},\beta))$,  the total sum of coefficients of negative powers of $z$ on the right-hand side must be zero (even though individual terms might have summands with negative powers of $z$). 
 
  The normal bundles of $F_{\pm}$ are of rank 1 and of weights $z$ and $-z$, we therefore obtain 
 \[\frac{1}{e_{\BC^*}(\BM^\vee_{-})}= -\frac{1}{z}+O(1/z^2), \quad \frac{1}{e_{\BC^*}(\BM_{+})}= \frac{1}{z}+O(1/z^2). \]
   By  Theorem \ref{Masterproper}, the master space is proper, hence after integrating the classes in the sum above, and taking the residue at $z=0$ of the resulting localisation formula (i.e., the coefficients of $1/z$),   we obtain the following relation,
\begin{multline*} 
\int_{[ \vec{A}dm^{\epsilon_-}_{g,N}(X^{(n)},\beta)]^{\mathrm{vir}}} \alpha -\int_{[ \vec{A}dm^{\epsilon_+}_{g,N}(X^{(n)},\beta)]^{\mathrm{vir}}} \alpha   
=    \mathrm{Res}_{z} \sum_{\underline{\beta}} \int_{[\vec{F}_{\underline{\beta}}]^{\mathrm{vir}}}  \frac{\alpha}{e_{\BC^*}(N_{\vec{F}_{\underline{\beta}}}^{\mathrm{vir}})},
\end{multline*}
where we used the description of the fixed loci from the beginning of Section \ref{secwall}. 
The analysis of the residue on the right-hand side presented in \cite[Section 7]{YZ} applies to our case, using Proposition \ref{isomophism}.  The resulting formula is the one claimed in the statement of the theorem.  \qed
\\

We define 
\[F^{\epsilon}_{g}(\mathbf{t}(z))=\sum^{\infty}_{N=0}\sum_{\beta}\frac{q^{\beta}}{N!}\langle \mathbf{t}(\psi), \dots, \mathbf{t}(\psi) \rangle^{\epsilon}_{g,N,\beta},\]
where $\mathbf{t}(z) \in H^{*}(\vec{\CI} X^{(n)},\BQ)[z]$ is a generic element, and the unstable terms are set to be zero. By repeatedly applying Theorem \ref{wallcrossingSym} crossing all walls from $0^+$ to $1$, we obtain the following result.
\begin{cor} \label{changeofvariable} For all $g\geq 1$ we have
\[F^{0^+}_{g}(\mathbf{t}(z))=F^{1}_{g}(\mathbf{t}(z)+I_+(-z)).\]
For $g=0$, the same equation holds modulo constant and linear terms in $\mathbf{t}(z)$. 	
\end{cor}
For $g = 0$, the relation is true only moduli linear terms in $\mathbf{t}(z)$, because the moduli space $Adm_{0,1}^{\epsilon_-}(X^{(n)}, \beta)$ is empty, if $ \deg(\beta) \leq 1/\epsilon_-$, and the wall-crossing takes a different form. See  \cite[Section 6.4]{YZ} for the statement in the case of quasimaps. 
\section{Del Pezzo} \label{del Pezzo}
\subsection{$I$-functions for del Pezzo surfaces} 
In this section, we determine the $I$-function for a del Pezzo surface $X=S$, and then evaluate the wall-crossing formula.

First, consider the expansion 
\[[zI(z)-z]_{+}=I_{1}(q)+(I_{0}(q)-1)z+I_{-1}(q)z^2+I_{-2}(q)z^3+\dots,\]
we will show that by the dimension constraint the terms $I_{-k}$ vanish for all $k\geq 1$. The virtual dimension of the space of relative maps with the ramification profile $\mu$, 
\[
\Mbar_{\mathsf{m}}(X\times \p^1/X_{\infty}, (\gamma,n),\mu) \subseteq \Mbar_{\mathsf{m}}(X\times \p^1/X_{\infty}, (\gamma,n)) 
\]  
is equal to 
\[\int_{\mathrm{c}_1(S)}\beta + n + \ell(\mu).\] 
Hence, by the virtual localisation, the classes involved in the definition of $I$-function,
\[\evsf_{*}\left( \frac{[\vec{F}_{\beta,\mu}]^{\mathrm{vir}}}{e_{\BC^{*}}(N^{\mathrm{vir}}_{\vec{F}_{\beta,\mu}})}\right)  \in H^*(S^{\mu},\BQ)[z^\pm],\] 
have the naive cohomological degree\footnote{That is the degree without the Chen--Ruan shifting discussed in (\ref{orbifoldcoh}).} equal to
\begin{equation} \label{quantity}
-2\left(\int_{\mathrm{c}_1(S)}\beta + n-\ell(\mu)\right).
\end{equation}
Since $S$ is a del Pezzo surface, the above quantity is non-positive, which implies that 
\[ I_{-k}(q)=0\] 
for all $k\geq 1$, because cohomology is non-negatively graded. Moreover, the quantity (\ref{quantity}) is zero, if and only if
\begin{itemize}
\item$ \mu=(1,\dots,1), \quad  \beta=(0,\mathsf m), \textrm{ or } $ 
\item $\mu=(1,\dots, 1), \quad \beta=(\gamma, \mathsf m), \textrm{ such that } \gamma \cdot \mathrm{c}_1(S)=1. $
\end{itemize}

In the latter case, $I_1(q)\in H^0(S^n)[\![q]\!]$, hence it will not contribute to the primary invariants by the string equation Lemma \ref{stringdivisor}, if $N>2$. We will therefore drop it and refer to the most recent version of  \cite{N} for the analysis of this part of $I$-function on the side of quasimaps. 

\subsection{The main component} From now on, we assume that 

\[ \mu=(1,\dots,1), \quad  \beta=(0,\mathsf m).  \]
Let us study $F_{\beta,\mu}$ for these values of $\mu$ and $\beta$. It is more convenient to put a standard order on ramification points over $\infty \in \p^1$, so let $\vec F_{\beta, \mu}$ be the resulting space. The full description of $\vec F_{\beta, \mu}$ is sketched in Section \ref{why}. We start with the simplest type of components of $\vec F_{\beta, \mu}$, which are the only ones that contribute to $I_{+,\beta}$. These components take the following form: \begin{equation} \label{iota}
\iota_i\colon \Mbar_{\mathsf h,p_i}\times S^n \hookrightarrow \vec  F_{\beta, \mu},
\end{equation}
where $\Mbar_{\mathsf h,p_i}$ is the moduli spaces of stable genus-$\mathsf h$ curve with \textit{one} marking labelled by $p_i$, $i=1,\dots N$. The embedding $\iota_i$ is constructed as follows.  Given a point
\[( (C,\mathbf x), x_1,\dots,x_n)) \in \Mbar_{\mathsf h,p_i}\times S^n,\] 
let 
\begin{equation} \label{rational}
(\tilde{P},p_1,\dots, p_n)=\coprod^{i=n}_{i=1} (\p^1,0)
\end{equation}
be an ordered disjoint union of $\p^1$ with markings at $0\in \p^1$. We define a curve $P$ by gluing $(\tilde{P},p_1,\dots, p_n)$ with $(C, p_i)$ at the marking with the same labelling. We define 
\[f_{\p^1}\colon P \rightarrow \p^1\] to be an identity on 
$\tilde P$ and contraction on $C$,  and
\[f_{S}\colon P \rightarrow S\] by mapping the $j$-th $\p^1$ in $P$ to the point $x_j \in S$; the curve $C$ is contracted together with the $i$-th $\p^1$ to $x_i$.  We thereby defined the inclusion
\[\iota_i((C,p), x_1,\dots,x_n))= (P, \p^1,f_{\p^1}\times f_S ),\] 
where the map $f_{\p^1}$ is clearly admissible at $\infty \in \p^1$.

By Lemma \ref{RHformula}, 
\begin{equation}\label{hm}
\mathsf h= \mathsf m/2,
\end{equation}
in particular, $\mathsf m$ is even. More generally, any connected component of $\vec F_{\beta, \mu}$ admits a similar description with the difference that there might be more markings on possibly disconnected $C$ by which it attaches to $\widetilde{P}$, i.e., $P$ has more nodes. These components of $\vec F_{\beta, \mu}$ are  irrelevant for the purpose of determining the truncation of $I$-function, we refer to Section \ref{why} for more details. 
 \subsection{Hodge integrals} Let us now consider the virtual fundamental classes and the normal bundles of the components $\Mbar_{\mathsf h,p_i}\times S^n$ considered above. 
The obstruction theory of contracted components is governed by the Hodge bundle $\BE$ on  $\Mbar_{\mathsf h,p_i}$, see \cite{FabP} for more details. In our case, the obstruction bundle of  $\Mbar_{\mathsf h,p_i}\times S^n$ is given by 
\[  (\pi^*_iT_S\otimes p^*\BE^{\vee}_{\mathsf h})\cdot \frac {e(\BE^{\vee}z)}{z(z-\psi_1)},\]
where 
$\pi_i\colon \Mbar_{\mathsf h,p_i}\times S^n \rightarrow S$ is the projection to $i$-th factor of $S^n$ and $p\colon \Mbar_{\mathsf h,p_i}\times S^n \rightarrow \Mbar_{\mathsf h,p_i} $ is the projection to $\Mbar_{\mathsf h,p_i}$; $e(V)$ denotes the Euler class of a vector bundle. We therefore have to determine the following classes 
\[\pi_*\left( e(\pi^*_i T_S\otimes p^*\BE^{\vee}_{\mathsf h})\cdot \frac {e(\BE^{\vee}z)}{z(z-\psi_1)}\right) \in H^{*}(S^n,\BQ)[z^{\pm}],\]
where $\pi \colon \Mbar_{h,p_i}\times S^n \rightarrow S^n$ is the natural projection, which is identified with evaluation morphism $\evsf$ via the inclusion (\ref{iota}).

Let $\ell_1$ and $\ell_2$ be the Chern roots of $\pi_i^*T_S$. Then we can rewrite the classes above as follows: 
\[\int_{\Mbar_{\mathsf h,1}}\frac{\Lambda^{\vee}(\ell_1)\cdot \Lambda^{\vee}(\ell_2) \cdot  \Lambda^{\vee}(z)}{z(z-\psi_1)},\]
where 
\[\Lambda^\vee(z):=e(\BE^{\vee}z)= \sum^{j=\mathsf h}_{j=0}(-1)^{\mathsf h-j}\lambda_{\mathsf h-j}z^j,\]
and similarly for $\Lambda^\vee(\ell_{1})$ and $\Lambda^\vee(\ell_{2})$.  

By putting these Hodge integrals into a generating series, we can explicitly determine them.  Note that below we sum over the degree $\mathsf m$ of the branching divisor, which in this case is related to the genus $\mathsf h$  by (\ref{hm}). The result was kindly communicated to the author by Maximilian Schimpf.
\begin{prop}[Maximilian Schimpf] \label{Max} We have
\[1+\sum_{\mathsf h>0} q^{2\mathsf h} \int_{\Mbar_{\mathsf h,1}}\frac{\Lambda^{\vee}(\ell_1)\cdot \Lambda^{\vee}(\ell_2) \cdot  \Lambda^{\vee}(z)}{z(z-\psi_1)}=\left( \frac{\mathrm{sin}(q/2)}{q/2}\right)^{\frac{\ell_1+\ell_2}{z}}.\]
\end{prop}
\textit{Proof.} The claim follows from the results of \cite{FabP}.  Firstly, by replacing every class $A$ that appears in the expression under the integral sign by $ z^{\deg_\BC(A)}A$ and multiplying the entire expression by $z^{-\dim(\Mbar_{\mathsf h,1})}$,  the integral remains unchanged; $\deg_\BC(A)$ denotes the cohomological degree of the class $A$. We therefore  obtain that
\begin{multline*}
1+\sum_{\mathsf h>0} q^{2\mathsf h} \int_{\Mbar_{\mathsf h,1}}\frac{\Lambda^{\vee}(\ell_1)\cdot \Lambda^{\vee}(\ell_2) \cdot  \BE^{\vee}(z)}{z(z-\psi_1)}\\
=1+\sum_{\mathsf h>0} q^{2\mathsf h} \int_{\Mbar_{\mathsf h,1}}\frac{\Lambda^{\vee}(\ell_1/z)\cdot \BE^{\vee}(\ell_2/z) \cdot  \BE^{\vee}(1)}{1-\psi_1}.
\end{multline*}
Now let 	\[a=\ell_1/z, \quad b=\ell_2/z\]and 
\[F(a,b)=1+\sum_{\mathsf h>0} q^{2\mathsf h} \int_{\Mbar_{\mathsf h,1}}\frac{\Lambda^{\vee}(a)\cdot \Lambda^{\vee}(b) \cdot  \Lambda^{\vee}(1)}{1-\psi_1}.\]By using virtual localisation on a moduli space of stable maps to $\p^1$, we obtain the following identities:
\begin{align*}
&F(a,b)\cdot F(-a,-b)=1, \\
&F(a,b)\cdot F(-a,1-b)=F(0,1).
\end{align*}
These identities, with the fact $F(a,b)$ is symmetric in $a$ and $b$, imply that 
\begin{equation} \label{identityf}
F(a,b)=F(0,1)^{a+b}
\end{equation}
for integer values of $a$ and $b$. Each coefficient of a power of $q$ in $F(a,b)$ is a polynomial in $a$ and $b$, hence the identity (\ref{identityf}) is in fact a functional identity. 

By the discussion in \cite[Section 2.2]{FP}  and by \cite[Proposition 2]{FP}, we obtain that 
\[F(0,1)= \frac{\mathrm{sin}(q/2)}{q/2},\]the claim now follows. 
\qed

\subsection{Other components} \label{why} Other components of $\vec F_{\beta, \mu}$ for $\mu=(1,\hdots,1)$ will consist  of  $\Mbar_{\mathsf h,k}$ with $k>0$, such that the curves are attached to $k$ disjoint copies of $\p^1$, or, more generally, products of such $\Mbar_{\mathsf h,k}$. In the former case, the obstruction bundle of contracted maps will be of the following form: 
\[
e(\pi^*T_S\otimes p^*\BE^{\vee}_{\mathsf h})\cdot  e(\BE^{\vee}z) \cdot    \frac{1}{z \prod^k_{i=1}(z-\psi_i)},
\]
where the rightmost factors arise due to extra nodes at which a curve is attached to $\p^1$ components (more precisely, from the obstruction to smooth them). As for $k=1$, $I$-functions will therefore be given by the following integral: 
\[\int_{\Mbar_{\mathsf h,k}}\frac{\Lambda^{\vee}(\ell_1)\cdot \Lambda^{\vee}(\ell_2) \cdot  \Lambda^{\vee}(z)}{z\prod^k_{i=1}(z-\psi_i)}.\] 
  By the dimension count, the biggest $z$-power of the integral above is $-k$. Hence for $k\geq 2$, they are not relevant for our computations. Recall that we are interested in $I_+(z)$, which is the non-polar part of $zI(z)-z$.  

 When there are several copies of $\Mbar_{\mathsf h,k}$, corresponding to several curves attached to $\p^1$ components, we will have a product of integrals above, which leads to the same conclusion.

\subsection{Evaluation of the wall-crossing formula} Overall, by taking the coefficient of $1/z$ in the formula from Proposition \ref{Max} and summing over  $n$ different choices of $p_i$ in  $\Mbar_{\mathsf h,p_i}$, we obtain
\begin{equation} \label{I1} I_1(q)=\log\left(\frac{\mathrm{sin}(q/2)}{q/2}\right)\cdot \sum^{n}_{i=1} (\mathbbm{1} \otimes \hdots  \mathrm{c}_1(S)  \hdots \otimes \mathbbm{1})+I_1(q)_{|\deg=0},
\end{equation}
where $I_1(q)_{|\deg=0}$ is the part of the $I$-function contained in $H^0(S^n)$, and we sum over all ways to put the class $\mathrm{c}_1(S)$ among identity classes in  $(\mathbbm{1}\otimes \dots   \mathrm{c}_1(S)  \dots \otimes \mathbbm{1})$, which correspond to the choice of the index of $p_i$.

We define
\[\left\langle \gamma_1, \ldots, \gamma_N \right\rangle^{ \epsilon}_{g, \gamma}(q):=\sum_{\mathsf{m}} \langle \gamma_1, \dots, \gamma_N \rangle_{g,(\gamma,\mathsf m)}^{\epsilon}q^{\mathsf m},\]
setting invariants corresponding to unstable values of $g$,$N$ and $\beta$ to zero.  By repeatedly applying Theorem \ref{wallcrossingSym}, we obtain that 
\[\left\langle \gamma_1, \ldots, \gamma_N \right\rangle^{ 0^+}_{g, \gamma}(q)=\sum_{k\geq 1}\frac{1}{k!}\left\langle \gamma_1, \ldots, \gamma_N, \underbrace{I_1(q),\dots, I_1(q)}_{k} \right\rangle^{1}_{g, \beta}(q).\]
To evaluate the expression above, we use the following divisor and string equations. 

\begin{lemma}  \label{stringdivisor} Assume $N >2$. Let $ \CD=\sum^{n}_{i=1} (\mathbbm{1}\otimes \hdots D \hdots \otimes \mathbbm{1}) \in H^*(X^n,\BQ) \subseteq H^*(\vec{\CI} X^{(n)},\BQ)$ be a divisor class, where the sum is taken over $n$ ways to put $D \in H^2(X,\BQ)$ among identity classes. Then the divisor equation holds, 
	\[ 
	\left\langle \gamma_1, \ldots, \gamma_N, \CD \right\rangle^{1}_{g, (\gamma,\mathsf{m})}=(\gamma \cdot D) \left\langle \gamma_1, \ldots, \gamma_N \right\rangle^{1}_{g, (\gamma,\mathsf{m})}.
	\]
	Let $\mathbbm{1} \in H^0(X^n)\subseteq H^*(\vec{\CI} X^{(n)},\BQ)$ the identity class, then the string equation holds, 
	\[ 
	\left\langle \gamma_1, \ldots, \gamma_N, \mathbbm{1} \right\rangle^{1}_{g, (\gamma,\mathsf{m})}=0
	\] 
	\end{lemma}
\textit{Proof.} Like in the case of the standard divisor equation, we use the forgetful map with respect to the marking on the target curve, 
\begin{equation} \label{for}
\vec{A}dm^{1}_{g,N+1}(X^{(n)},\beta) \rightarrow \vec{A}dm^{1}_{g,N}(X^{(n)},\beta),
\end{equation}
which exists by the assumption on $N$ (if $N=2$ and $\beta=0$, then the forgetful map does not exist due to the instability of curves). 
The  fiber of the forgetful map over $[(P,C, \mathbf{x},f)] \in \vec{A}dm^{1}_{g,N}(X^{(n)},\beta)$  is the degree $n!$ ramified cover $C$ associated to the degree $n$ cover $P\rightarrow C$ by taking the $S^n$-torsor away from ramification points (given by ordering the fibers) and applying the Riemann existence theorem over the ramification points (in other words, it is the $S_n$-torsor associated to the degree $n$ \'etale cover on the corresponding twisted curve). The degree of this $n!$ cover inside $X^n$ is  \[(n-1)! \sum^{i=n}_{i=1} (\mathrm{pt}\otimes \hdots   \gamma  \hdots \otimes \mathrm{pt}),\]
since $P \rightarrow C$ is a $S_{n-1}$-quotient of  the degree $n!$ cover, and the class of $P$ inside $X$ is $\gamma$ (see also Section \ref{Relation1} for more details on how degree $n$  \'etale covers and $S_n$-torsors are related). 
In particular, we note that 
\[  (n-1)! \sum^{n}_{i=1} (\mathrm{pt}\otimes \hdots   \gamma \hdots \otimes \mathrm{pt}) \cdot \sum^{n}_{i=1} (\mathbbm{1}\otimes \hdots D  \hdots \otimes \mathbbm{1})  =n! (\gamma\cdot D), \]
and $|\Aut((1,\dots, 1))|=n!$. 
Hence by the projection formula applied to the forgetful map (\ref{for}), we obtain the divisor equation. The string equation is proved in the same way. 
\qed  

\begin{rmk} Note that, in general, the divisor equation does not hold for the orbifold GW theory. For Lemma \ref{stringdivisor}, it is crucial to assume that the divisor $\CD$ is contained in the untwisted sector, i.e., $\mu=(1,\dots, 1)$. Indeed, otherwise, fibers of the forgetful map do not admit the same description.  The failure of the divisor equation for $\mu=(2,1,\hdots,1)$ gives rise to the change of variables $y=-e^{iq}$ in the crepant resolution conjecture, as explained in Section \ref{theremark}.  
	\end{rmk}

 Using Lemma  \ref{I1}, Theorem \ref{wallcrossingSym} gives us the following result.  We make the assumption $N>2$, so that the divisor and string equations can be applied in degree $\beta=0$.  
\begin{cor} \label{change} Assuming $N>2$ and $S$ is a del Pezzo surface, we have
\[\left\langle \gamma_1, \ldots, \gamma_N \right\rangle^{0^+}_{g, \gamma}(q)= \left( \frac{\mathrm{sin}(q/2)}{q/2} \right)^{\gamma  \cdot \mathrm{c}_1(S)} \cdot  \left\langle \gamma_1, \ldots, \gamma_N \right\rangle^{1}_{g, \gamma}(q).\]
\end{cor}

 We also note that Corollary \ref{equivalent} in fact holds for all surfaces, if $\gamma=0$, as the calculation of the $I$-function applies to classes $(0, \mathsf{m})$ verbatim. 
 
 \begin{cor} \label{change2} Assuming $N> 2$ and $S$ is a projective surface, we have
 	\[\left\langle \gamma_1, \ldots, \gamma_N \right\rangle^{0^+}_{g, \gamma}(q)=  \left\langle \gamma_1, \ldots, \gamma_N \right\rangle^{1}_{g, \gamma}(q).\]
 \end{cor}

\section{Crepant resolution conjecture} \label{crepant}

\subsection{Changing $\epsilon$-admissibility} In what follows, there will be two types of $\epsilon$-stabilities, the one associated to quasimaps from \cite{N}, and the one considered in this work. In both case, $\epsilon$ lives in $\BR_{>0}$. To distinguish the two, we let $\epsilon$ live in $[-1,0) \in \BR_{<0}$ in the case of $\epsilon$-admissibility. Concretely, this amounts to introducing negative signs in inequalities in Definition \ref{epsilonadm}. In particular, the side of  admissible covers to $X\times C$ will correspond to $\epsilon=-1$ now, while the side of stable maps to $X\times C$ will correspond to $\epsilon=0^-$. Positive $\epsilon$ will be reserved for quasimaps from now on, 
\begin{align*}
	\epsilon \text{-admissibility}\in&   \BR_{<0}  \\
	 \epsilon\text{-stability of quasimaps} \in& \BR_{>0}.  
	\end{align*}
\begin{rmk} A reader might wonder why $\epsilon$-stability for quasimaps lives in the entire $\BR_{>0}$, while $\epsilon$-admissibility only in $(0,1]$. In fact, the last wall for $\epsilon$-stability of  quasimaps is at $\epsilon=1$, and for all $\epsilon >1$, $\epsilon$-stable quasimaps are just stable quasimaps, (at least for $(g,N)\neq(0,1)$). It is however customary to  use the entire $\BR_{>0}$ nevertheless. The fact that we do not cross the wall at $\epsilon=1$ for $\epsilon$-admissibility (i.e., we do not forbid simple ramifications) is a feature of $\epsilon$-admissibility and maps to symmetric products.  The number of simple ramifications constitutes the extended degree of maps \cite{BG}. 
	\end{rmk}
\subsection{Classes on Hilbert schemes} To a cohomology-weighted partition 
\[\vec{\mu}=((\mu_1,\delta_{\ell_1}), \dots, (\mu_k,\delta_{\ell_k}))\]we can also associate a class in $H^*(S^{[n]},\BQ)$, using Nakajima operators, 
\[\theta(\vec{\mu}):=\frac{1}{\mathfrak{z}(\mu)}P_{\delta_{\ell_1}}[\mu_1]\dotsb P_{\delta_{\ell_k}}[\mu_k]\cdot 1_S \in H^*(S^{[n]},\BQ),\]
where operators are ordered according to the standard ordering (see Subsection \ref{Classes}). For more details on these classes, we refer to \cite[Chapter 8]{N99}, or to \cite[Section 0.2]{Ob18} in a context more relevant to us. 

Classes $\theta(\vec{\mu})$ and $\lambda(\vec{\mu})$ provide the basis for the corresponding cohomologies. Moreover, both $H^*(S^{[n]},\BC)$ and $H_{\mathrm{orb}}^{*}(S^{(n)},\BC)$ are equipped with non-degenerate pairings, such that

\begin{equation}\label{pairing}
\langle \theta(\vec{\mu}), \theta(\vec{\eta}^\vee) \rangle_{S^{[n]}}=\frac{(-1)^{\mathrm{age}(\mu)}\delta_{\mu,\eta}}{\mathfrak{z}(\mu)}, \quad \langle \lambda(\vec{\mu}), \lambda(\vec{\eta}^\vee) \rangle_{S^{(n)}}=\frac{\delta_{\mu,\eta}}{\mathfrak{z}(\mu)},
\end{equation}
where $\vec{\eta}^\vee$ is the dual cohomology-weighted partition - the classes $\delta_i$ are replaced with their cohomological duals $\delta_i^\vee$, while $\eta$ stays the same; and  $\delta_{\mu,\eta}$ is the Kronecker delta function.

\begin{prop}\label{HilbSym} There exists an isomorphism of graded vector spaces which respects the pairing,
\[L\colon H_{\mathrm{orb}}^{*}(S^{(n)},\BC)\simeq H^*(S^{[n]},\BC),\]
\[	L(\lambda(\vec{\mu}))= (-i)^{\mathrm{age}(\mu)}\theta(\vec{\mu}), \]
such that the grading on the left hand-side is given by (\ref{orbifoldcoh}).
\end{prop}
\textit{Proof.} See \cite[Proposition 3.5]{FG}.
\qed 
\begin{rmk} The peculiar choice of the identification with a factor $(-i)^{\mathrm{age}(\mu)}$ is justified by the crepant resolution conjecture - it makes the above isomorphism respect the pairings on both sides.  See also \cite[Section 3]{Che}.
	
	The square of this factor is also used to modify the product in \cite[Definition 3.9]{FG}. Note that in their notation, $\ell(g)$ of an element $g \in S_n$ means $\mathrm{age}(g)$ in our notation (after associating a partition $\mu$ to the conjugacy class of $g$).

\end{rmk}
\subsection{Quasimaps and admissible covers} \label{qausiadm}  Let $S^{[n]}$ be the Hilbert scheme of $n$-points on $S$. In this section, we use the quasimap theory to $S^{[n]}$ of \cite{N}. The Nakajima basis provides the following identification,
\begin{equation}\label{degreechern} H_2(S^{[n]},\BZ)\cong H_2(S,\BZ)\oplus \BZ\cdot A,
\end{equation}
where $A:=\mathfrak{z}(\mu)\theta(\vec{\mu})$ for $\vec{\mu}=((2,\mathrm{pt}), (1,\mathrm{pt}),\dots (1,\mathrm{pt}))$. 
We refer to \cite[Section 8.2]{N} for how to associate a class in $H_2(S^{[n]},\BZ)$ to a degree of quasimap, using the Nakajima basis.\footnote{Note that the classes with respect to the Nakajima basis are denoted by $m$ in \cite[Section 8.2]{N}.}
Given classes $\gamma_i \in H^*_{\mathrm{orb}}(S^{(n)},\BC)$, $i=1,\dots N$, and a class \[(\gamma,\mathsf k) \in  H_2(S,\BZ)\oplus \BZ,\]
such that the integer $\mathsf{k}$ is  the degree of quasimaps with respect to the exceptional curve class $A$. For $\epsilon \in \BR_{>0}\cup \{\infty\}$ we set
\[\langle \gamma_1, \dots, \gamma_N \rangle^{\epsilon}_{g,(\gamma,\mathsf{k})}:= \langle L(\gamma_1), \dots, L(\gamma_N) \rangle^{\sharp,\epsilon}_{g,(\gamma,\mathsf{k})} \in \BC,\]
the invariants on the right are defined in \cite[Definition 6.9]{N} and $L$ is defined in Proposition \ref{HilbSym}. The integer $\mathsf{k}$ is the degree of a quasimap with respect to the exceptional curve class. We set
\[\left\langle \gamma_1, \ldots, \gamma_N \right\rangle^{ \epsilon}_{g, \gamma}(y):=\sum_{\mathsf k}  \langle \gamma_1, \dots, \gamma_N \rangle_{g,(\gamma,\mathsf{k})}^{\sharp,\epsilon}y^\mathsf{k}.\]
For $\epsilon=0^+$, these are the relative PT invariants of the relative geometry $S\times C_{g,N}\rightarrow \Mbar_{g,N}$. The summation over $\mathsf{k}$ with respect to the identification (\ref{degreechern}) corresponds to the summation over $\ch_3+\mathrm{c}_1(S)\cdot \gamma/2$ of a subscheme.
\\

Using wall-crossings, we will now show the compatibility of \textsf{PT/GW} and  \textsf{CRC}.  Let us firstly spell out our conventions.

\begin{itemize}
\item We sum over the degree of the branching divisor instead of the genus of the source curve. Assuming $\gamma_i$'s are homogenous with respect to the decomposition  $\vec{\CI}S^{(n)}:=\coprod_{\mu} S^{\mu} $, the genus $\mathsf{h}$ and the degree $\mathsf{m}$ are related by Lemma \ref{RHformula},
\[2\mathsf{h}-2=n(2g-2)+\mathsf m + \sum \mathrm{age}(\gamma_i),\]
where $\mathrm{age}(\gamma_i)$ is equal to $\mathrm{age}(\mu)$ for a partition $\mu$ such that $\gamma_i \in S^{\mu}$. 
For $\epsilon \in [-1,0)$, let
\['\left\langle \gamma_1, \dots, \gamma_N \right\rangle^{\epsilon}_{g, \gamma}(q):= \sum_{\mathsf h}  \langle \gamma_1, \dots, \gamma_N \rangle_{g,(\gamma,\mathsf{h})}^{\epsilon}q^\mathsf{2h-2}\]
be generating series where the summation is taken over the genus instead. Then the two generating series are  related as follows:
\['\left\langle \gamma_1, \dots, \gamma_N \right\rangle^{ \epsilon}_{g, \gamma}(q) = q^{n(2g-2)+\sum \mathrm{age}(\gamma_i)} \cdot \left\langle \gamma_1, \dots, \gamma_N \right\rangle^{ \epsilon}_{g, \gamma}(q).\]
\item We sum over the degree $\mathsf{k}$ with respect to the exceptional divisor on $S^{[n]}$  instead of Euler characteristics $\chi$. For $\epsilon \in \BR_{>0 }\cup \{ \infty\}$, let
\['\left\langle \gamma_1, \dots, \gamma_N \right\rangle^{\epsilon}_{g, \gamma}(y):= \sum_{\mathsf \chi}  \langle \gamma_1, \dots, \gamma_N \rangle_{g,(\gamma, \chi)}^{\sharp,\epsilon}y^\chi\]
be the generating series where the summation is taken over Euler characteristics $\chi$. Then by \cite[Section 8.2]{N}, the two generating series are related as follows
\['\left\langle \gamma_1, \dots, \gamma_N \right\rangle^{\epsilon}_{g, \gamma}(y) = y^{n(1-g)} \cdot \left\langle \gamma_1, \dots, \gamma_N \right\rangle^{\epsilon}_{g, \gamma}(y).\]
\item The identification of Proposition \ref{HilbSym} has a factor of $(-i)^{\mathrm{age}(\mu)}$. 
\end{itemize}
Taking into account all the conventions above,  the PT/GW correspondence, \cite[Conjectures 2R, 3R]{MNOP}, can be expressed\footnote{We take the liberty to extend the statement of the conjecture in \cite{MNOP} from a fixed curve to a moving one, and from one relative insertion to multiple ones; see \cite{PP} for the statement of the conjectures with multiple relative insertions.} as follows. Recall that $\mathrm{c}_1(S\times C)=\mathrm{c}_1(S)+(2-2g)$, and the degree of a curve $P$ over $S\times C$ is $(\gamma,n)$. 
\\

\noindent
\textbf{PT/GW.} \textit{The generating series} $\left\langle \gamma_1, \dots, \gamma_N \right\rangle^{0^+}_{g, \gamma}(y)$ \textit{is a Taylor expansion of a rational function around $0$},\textit{ such that under the change of variables }$y=-e^{iq}$,
\[(-y)^{-\gamma\cdot \mathrm c_1(S)/2}\cdot \left\langle \gamma_1, \dots, \gamma_N \right\rangle^{0^+}_{g, \gamma}(y)=(-iu)^{\gamma\cdot \mathrm c_1(S)} \cdot \left\langle \gamma_1, \dots, \gamma_N \right\rangle^{0^-}_{g, \gamma}(q).\]

Note that the additional factor $(-i)^{-\mathrm{age}(\mu)}$ which appears on the left hand side of \cite[Conjecture 3R]{MNOP} was absorbed by the same factor in the definition of the isomorphism $L$.  

On the other hand, the crepant resolution conjecture of \cite{BG} takes the following form. 
\\

\noindent\textbf{CRC.} \textit{The generating series} $\left\langle \gamma_1, \dots, \gamma_N \right\rangle^{\infty}_{g, \gamma}(y)$\textit{ is a Taylor expansion of a rational function around 0, such that under the change of variables} $y=-e^{iq}$,
\[\left\langle \gamma_1, \dots, \gamma_N \right\rangle^{\infty}_{g, \gamma}(y) =\left\langle \gamma_1, \dots, \gamma_N \right\rangle^{-1}_{g, \gamma}(q).\]

Assume now that $S$ is a del Pezzo surface. Let us apply our wall-crossing formulas. Using Corollary \ref{change}, we obtain
\begin{equation*}
(-iu)^{\mathrm{c}_1(S)\cdot \gamma} \cdot \left\langle \gamma_1, \dots , \gamma_N \right\rangle^{0^-}_{g, \gamma}(q)=(e^{-iq/2}-e^{iq/2})^{\gamma \cdot \mathrm c_1(S)}\cdot \left\langle \gamma_1, \dots, \gamma_N \right\rangle^{-1}_{g, \gamma}(q).
\end{equation*}
Using \cite[Corollary 7.7]{N}, we obtain 
\begin{multline*}
(-y)^{-\gamma\cdot \mathrm c_1(S)/2} \left\langle \gamma_1, \dots, \gamma_N \right\rangle^{0^+}_{g, \gamma}(y)=\\
(-y)^{-\gamma\cdot \mathrm c_1(S)/2}\cdot (1+y)^{\gamma \cdot \mathrm c_1(S)}\cdot \left\langle \gamma_1, \dots, \gamma_N \right\rangle^{\infty}_{g, \gamma}(y).
\end{multline*}

Substituting $y=-e^{iq}$ into the above equations, we find that by both wall-crossings the statements of \textsf{PT/GW} and \textsf{CRC} in the form presented above are equivalent. 
\begin{cor} \label{equivalent} Assuming $N>2$ and $S$ is a del Pezzo surface, we have an equivalence of  the conjectures, 
\[\mathbf{PT/GW} \iff \mathbf{CRC}.\]
\end{cor}

\subsection{The change of variables} \label{theremark} In fact, the change of variables $y=-e^{iq}$ appears naturally through a weaker version of the crepant resolution conjecture. Firstly consider the following cohomology-weighted partition, 
	\[ (2,1^{n-2})=((2,\mathbb{1}), (1,\mathbb{1}), \dots, (1,\mathbb{1}))\ \]

The associated class $\theta(2,1^{n-1})$ on $S^{[n]}$ is the multiple of the exceptional divisor class\footnote{The factor of $2$  is due to the factor $\prod_j \mu_j$. } $E/2$.  Recall that $E$ is the locus of non-reduced $n$-points on $S$, which is also the exceptional divisor of the resolution of singularities $S^{[n]} \rightarrow S^{(n)}$. By (\ref{pairing}), the intersection number of the exceptional divisor $E/2$ and the exceptional curve class $A$ is $-1$. In particular, by the divisor equation, we have 
\[\left\langle \gamma_1, \dots, \gamma_N,\underbrace{E/2,\dots, E/2}_{b} \right \rangle^{\infty}_{g, \gamma,\mathsf{k}}=(-\mathsf{k})^b\left\langle \gamma_1, \dots, \gamma_N \right \rangle^{\infty}_{g, \gamma,\mathsf{k}}.\]
On the other side, inserting the class $I_{(2)}:=\lambda (2,1^{n-2})$ on $S^{(n)}$ has a very particular role - by definition, it just introduces a simple ramification, thereby raising $\mathsf{m}$. In particular, we have 

\[ \left\langle \gamma_1, \dots, \gamma_N, \underbrace{I_{(2)},\dots, I_{(2)}}_{b} \right \rangle^{-1}_{g, \gamma,\mathsf m}/\mathsf{m}!=\left\langle \gamma_1, \dots, \gamma_N \right \rangle^{-1}_{g, \gamma,\mathsf{m}+b}, \]
where we divide by the factorial, because the simple ramifications are not ordered. 
This discrepancy between classes on both sides is the reason why $y=-e^{iq}$ appears, as we now show.   

Let us sum over multiples of the exceptional curve class on the side of $S^{[n]}$ and do not sum over the degree of the branching divisor on the side of $S^{(n)}$, and assume that we have 
\begin{equation} \label{weak}
	\left\langle \gamma_1, \dots, \gamma_N \right\rangle^{\infty}_{g, \gamma}(-1) =\left\langle \gamma_1, \dots, \gamma_N \right\rangle^{-1}_{g, \gamma,0} 
	\end{equation}
for all insertions, using the isomorphism $L$. Then it must be true that under the specialisation $y=-1$ we have
\begin{multline*}
\sum_{k,\mathsf{m}}\left\langle \gamma_1, \dots, \gamma_N, E/2,\dots, E/2 \right \rangle^{\infty}_{g, \gamma,\mathsf{k}}y^\mathsf{k} (iq)^{\mathsf{m}}/\mathsf{m}! \mid_{y=-1}\\
=\sum_{k,\mathsf{m}}\left\langle \gamma_1, \dots, \gamma_N, I_{(2)},\dots, I_{(2)}\right \rangle^{-1}_{g, \gamma, 0}q^{\mathsf{m}}/\mathsf{m}!,
\end{multline*}
where $\mathsf{m}$ is the number of insertions corresponding to $E$ and $I_{(2)}$. 
Using the properties of these classes stated above, we obtain that 
\begin{equation}
	\left\langle \gamma_1, \dots, \gamma_N \right\rangle^{\infty}_{g, \gamma}(-e^{iq}) =\left\langle \gamma_1, \dots, \gamma_N \right\rangle^{-1}_{g, \gamma}(q).
\end{equation}
We therefore conclude that the weaker version of the crepant resolution conjecture  (\ref{weak}) is sufficient to obtain the version involving the change of variables $y=-e^{iq}$.

\subsection{Quantum cohomology} \label{qcoh} Let $g=0, N=3$. This is a particularly nice case, firstly because these invariants collectively are known as \textit{quantum cohomology}. Secondly, the moduli space of genus-0 curves with 3 markings is a point. Hence the invariants $\left\langle \gamma_1, \gamma_2, \gamma_3 \right\rangle^{0^+}_{0, \gamma}(y)$ are relative PT invariants of $S\times \p^1$ relative to the vertical divisor $S_{0,1,\infty}$, while $\left\langle \gamma_1, \gamma_2, \gamma_3 \right\rangle^{0^-}_{0, \gamma}(q)$ are relative GW invariants of $S\times \p^1$. 
In \cite{PP},  \textsf{PT/GW} is established for $S\times \p^1$ relative to $S\times\{0,1,\infty\}$, if $S$ is toric. Corollary \ref{equivalent} then implies the following.

\begin{cor} \label{quantum} If $S$ is toric del Pezzo, $g=0$ and $N=3$, then 
$\mathbf{CRC}$ holds in all classes. 
\end{cor}

The 3-point invariants together with a pairing define the quantum product $-\ast_q -$ on cohomologies,
\[\langle \gamma_1\ast_q \gamma_2,\gamma_3 \rangle := \sum_\beta \langle \gamma_1,\gamma_2,\gamma_3 \rangle_{0,\beta} q^\beta,\]
where $\beta \in H_2(S,\BZ)\oplus \BZ$. 
 By the above corollary and Proposition \ref{HilbSym}, we obtain that the identification $L$ respects the quantum ring structure on both sides after the change of variables $y=-e^{iq}$. 
 
We also note that for classes of the form $(0,\mathsf m)$, we recover the result of \cite{Che} for toric surfaces by Corollary \ref{change}.


\bibliographystyle{amsalpha}
\bibliography{QMs}
\end{document}